\documentclass[12pt]{amsart}

\usepackage{amsmath}
\usepackage{amsfonts}
\usepackage{amsthm}
\usepackage{amstext}
\usepackage{amssymb}
\usepackage{fullpage}
\usepackage{framed}
\usepackage{fancybox}
\usepackage{color}
\usepackage{mdwlist}
\usepackage{pifont}
\usepackage{graphicx}               
\usepackage{dsfont}
\usepackage{caption}
\usepackage{subcaption}
\usepackage{verbatim}
\usepackage{tikz}
\usepackage{enumerate}
\usepackage{cite}
\usepackage{newclude}
\usepackage{hyperref}
\usepackage{xcolor}
\hypersetup{
    colorlinks,
    linkcolor={red!70!black},
    citecolor={green!50!black},
    urlcolor={blue!80!black}
}

\usetikzlibrary{decorations.pathreplacing}
\usetikzlibrary{calc}
\usetikzlibrary{patterns}

\newcommand{\R}{\mathbb{R}}
\newcommand{\Z}{\mathbb{Z}}
\newcommand{\N}{\mathbb{N}}

\newcommand{\slr}{{\rm SL}_2(\R)}
\newcommand{\slnr}{{\rm SL}_n(\R)}
\newcommand{\slz}{{\rm SL}_2(\Z)}
\newcommand{\slnz}{{\rm SL}_n(\Z)}

\newcommand{\CX}{C^\infty_c(X)}
\newcommand{\sob}[2]{\mathcal{S}_{#1,#2}}
\newcommand{\sobt}[1]{\mathcal{S}_{#1}}
\newcommand{\Ismth}{I_\text{smth}}

\newcommand{\nn}[1]{\left|\!\left| #1 \right|\!\right|}

\newcommand{\dst}[1]{\displaystyle{#1}}

\newcommand{\bigzero}{\mbox{\normalfont\Large 0}}

\newcommand{\bigast}{\mbox{\normalfont\Large $\ast$}}

\theoremstyle{theorem}
\newtheorem{thm}{Theorem}[section]
\newtheorem{lem}{Lemma}[section]

\newtheorem{cor}{Corollary}[section]
\newtheorem*{cor*}{Corollary}
\newtheorem*{lem*}{Lemma}
\theoremstyle{definition}
\newtheorem{defn}{Definition}[section]
\newtheorem*{rem}{Remark}
\newtheorem*{rmk}{Remark}

\newtheorem{innercustomlem}{Lemma}
\newenvironment{customlem}[1]
  {\renewcommand\theinnercustomlem{#1}\innercustomlem}
  {\endinnercustomlem}

\title{Almost-Primes in Horospherical Flows on the Space of Lattices}
\author{
Taylor McAdam \\
Department of Mathematics \\
University of California, San Diego
}

\begin{document} 

\maketitle

\begin{abstract}
We study the asymptotic distribution of almost-prime entries of abelian horospherical flows on $\Gamma\backslash\slnr$, where $\Gamma$ is either $\slnz$ or a cocompact lattice.  In the cocompact case, we obtain a result that implies density of almost-primes of a sufficient order, and in the space of lattices we show the density of almost-primes in the orbits of points satisfying a certain Diophantine condition.  Along the way we give an effective equidistribution result for arbitrary horospherical flows on the space of lattices, as well as an effective rate for the equidistribution of arithmetic sequences of times in abelian horospherical flows.
\end{abstract}


\section{Introduction}

There is an intimate connection between number theory and dynamics on homogeneous spaces.  The case of $\Gamma\backslash\slnr$ where $\Gamma$ is a lattice is one particularly interesting and well-studied example, and when $\Gamma = \slnz$ this space can be identified with the space of unimodular lattices in $\R^n$.  It is well known, for instance, that the geodesic flow on $\slz\backslash\slr$ is related to Diophantine approximation of real numbers by rationals, which can be generalized to metric Diophantine approximation on manifolds (see \cite{KMNondiv}). Another famous example is Margulis's proof of the Oppenheim conjecture (see \cite{Opp}, \cite{MOpp}), which uses Ranghunathan's insight that the conjecture can be reduced to a statement about unipotent orbits in ${\rm  SL}_3(\Z)\backslash{\rm  SL}_3(\R)$.  Quantitative proofs of the Oppenheim conjecture are given in \cite{QuantOpp} and \cite{EMMQuantOpp}, and these make use of Ratner's work on measure rigidity for unipotent orbits (see \cite{Rat3}).
In addition to the contributions of dynamics to the field of number theory, there are also many number-theoretic questions that are of independent interest in dynamical systems, such as understanding the distribution of certain discrete subsets of times in a dynamical system (e.g. polynomial sequences or primes).  Many of these questions remain open---see, for example, the conjecture of Shah in the introduction of \cite{ShahConj} or the collection of conjectures by Margulis listed under Question 16 of \cite{Conjectures}.

Equidistribution results play an important role in dynamical systems and their applications to number theory.  Roughly speaking, a subset of some orbit is said to equidistribute with respect to a given probability measure if it spends the expected amount of time in subsets, i.e., if the proportion of the orbit landing within any set is given by the measure of that set.  The dual of this notion is that averages of any suitably nice function over larger and larger pieces of the orbit coverge weakly to the average of that function over the whole space with respect to the given measure.  Often in applications to number theory it is important that an equidistribution result be effective---that is, that there is a known rate of convergence.  Another question that can be asked is whether we can leverage the known equidistribution of a full orbit to obtain information about the distribution of certain ``sparse" subsets of that orbit.  Examples of research on sparse equidistribution problems can be found in \cite{VenkSparse}, \cite{SarnackUbis}, \cite{ShahConj}, and \cite{Heegner1} and \cite{Heegner2}.

Horospherical flows are a type of dynamical system arising naturally in the study of homogeous spaces.  A subgroup of a Lie group $G$ is said to be horospherical if it is contracted (conversely, expanded) under iteration of the adjoint action of some element of $G$ (see Section \ref{sect:horo} for a more precise definition).  It can be shown that any horospherical subgroup is unipotent, although not every unipotent subgroup can be realized as the horospherical subgroup corresponding to an element of $G$.  In general, horospherical flows are easier to study than more general unipotent flows, as the expansion property can be used along with dynamical information about the corresponding one-parameter subgroup to great effect.

Actions by horospherical and unipotent subgroups have been studied extensively.  
It was proved in \cite{Hedl} that the horocyclic flow on $\Gamma\backslash\slr$ for $\Gamma$ cocompact is minimal and later shown in \cite{Furst} to be uniquely ergodic.  These results were extended in \cite{Veech} and \cite{REWP} to more general horospherical flows on compact quotients of suitable Lie groups.  For $\Gamma$ non-uniform, we do not have unique ergodicity or minimality, however it was proved in \cite{Marg71} that orbits of unipotent flows cannot diverge to infinity in noncompact settings, which was refined in Dani's nondivergence theorem in \cite{Dani84a}, \cite{DaniNondiv}.  Moreover, it was shown in \cite{Dani78a} (for the case of $\Gamma\backslash\slr$ noncompact) and \cite{Dani81}, \cite{Dani86} (for more general noncompact homogeneous spaces) that horocyclic/horospherical flows have nice (finite volume, homogeneous submanifold) orbit closures and that every ergodic probability measure invariant under such a flow is the natural Lebesgue measure on some such orbit closure.
This work paved the way for a series of breakthrough papers, culminating in \cite{Rat3} and summarized in \cite{RatSummary}, in which Ratner resolved conjectures of Raghunathan and Dani by giving an essentially complete description of unipotent orbit closures and unipotent-invariant measures on homogeneous spaces.

More recently, many important results for horospheres and related actions have been effectivized.  Quantitative versions of Dani's nondivergence theorem were given in \cite{DaniMarg} and \cite{KMNondiv}, as well as a discrete version for $\slz\backslash\slr$ in \cite{SarnackUbis}.  In \cite{Burger}, Burger gave an effective rate for the equidistribution of the horocycle flow on compact quotients of $\slr$, which was improved upon and extended to the noncompact setting in \cite{FF} and \cite{StromDev}.  Many authors have also considered the effective equidistribution of closed horocyclic and horospherical orbits in a variety of settings, such as \cite{SarnackAsymp}, \cite{StromClosed}, \cite{KMExpanding}, \cite{LeeOh}, and \cite{KDabb}, although this list is by no means complete.  In studying both closed horospherical orbits and long pieces of generic horospherical orbits, one can make use of the ``thickening" argument developed by Margulis in his thesis \cite{Mthesis}.  This uses a known rate of mixing for the semisimple flow with respect to which the given subgroup is horospherical along with the expansion property to get a rate for the horospherical flow.  The key exponential rate for semisimple flows (and much more general actions) is given in \cite{KMBddOrbits}.  Other effective results of interest include \cite{EMV}, \cite{GreenTao1},  \cite{StromRat}, and \cite{EMMV}.  We note that most of these results use in some way a spectral gap for the action by translations of the ambient group $G$ or certain subgroups of $G$ on $L^2(\Gamma\backslash G)$.

One reason for wanting effective results is that many applications involving number theory require that the error in relevant approximations be controlled in a quantitative way.  For example, \cite{VenkSparse} makes use of effective equidistribution for the horocycle flow to 
derive an effective rate for the equidistribution of arithmetic sequences, which he then uses to show that sequences of integer times raised to small powers also equidistribute in $\Gamma\backslash\slr$ for $\Gamma$ cocompact.  As another example, \cite{SarnackUbis} uses effective equidistribution results along with sieving to demonstrate that prime times in the horocycle flow on $\slz\backslash\slr$ are dense in a set of positive measure.  More generally, if we hope to apply sieve methods to any equidistribution problem, we will need to have some way of quantitatively controlling the error.

In this paper, we are interested in the asymptotic distribution of almost-primes (i.e. integers having fewer than a fixed number of prime factors) in horospherical flows on the space of lattices and on compact quotients of $\slnr$.
Our main results are summarized in the following two theorems:

\begin{thm}\label{thm:CmptThm}
Let $\Gamma<\slnr$ be a cocompact lattice and  $u({\bf t})$ be an abelian horospherical flow on $\Gamma\backslash\slnr$ of dimension $d$.
Then there exists a constant $M$ (depending only on $n$, $d$, and $\Gamma$) such that for any $x\in\Gamma\backslash\slnr$, the set 
\[
\{xu(k_1, k_2, \cdots, k_d) \hspace{2pt}\vert\hspace{2pt} k_i\in\Z \text{ has fewer than } M \text{ prime factors}\}
\]
is dense in $\Gamma\backslash\slnr$.
\end{thm}

We remark that the dependence of the constant $M$ on $\Gamma$ arises from the spectral gap, so this dependence can be removed if $n\geq 3$ or for $n=2$ if $\Gamma$ is a congruence lattice.

For a horospherical flow $u({\bf t})$ on $\slnz\backslash\slnr$, we say that $x=\slnz g$ is \textit{strongly polynomially} $\delta$\textit{-Diophantine} if there exists some sequence $T_i\to\infty$ as $i\to\infty$ such that 
\[
\inf_{\substack{w\in \Lambda^j(\Z^n)\setminus\{0\}\\j=1, \cdots, n-1}} \sup_{{\bf t}\in [0,T_i]^d}\nn{wgu({\bf t})}>T_i^\delta
\]
for all $i\in\N$.

\begin{thm}\label{thm:NoncmptThm}
Let $u({\bf t})$ be an abelian horospherical flow on $\slnz\backslash\slnr$ of dimension $d$ and let $x\in\slnz\backslash\slnr$ be strongly polynomially $\delta$-Diophantine for some $\delta>0$.  Then there exists a constant $M_\delta$ (depending on $\delta$, $n$, and $d$) such that 
\[
\{xu(k_1, k_2, \cdots, k_d) \hspace{2pt}\vert\hspace{2pt} k_i\in\Z \text{ has fewer than } M_\delta \text{ prime factors}\}
\]
is dense in $\slnz\backslash\slnr$.
\end{thm}

A brief outline of the paper is as follows:  

In Section \ref{sect:prelims}, we establish the basic notation that will be used throughout the paper and introduce the key facts and theorems that we use in our analysis.  We also prove a small corollary of the nondivergence theorem in \cite{KMNondiv} that applies to the particular setting of this paper.

In Section \ref{sect:Equidistribution}, we prove an effective equidistribution result for long orbits of arbitrary horospherical flows on the space of lattices.  The proof makes use of the ``thickening" argument of Margulis, leveraging the exponential mixing properties of the subgroup with respect to which the flow of interest is horospherical, which is itself a consequence of a spectral gap.  The main result in this section is probably not surprising to experts, but the author was unable to locate a result in the literature that is stated in the way presented here.

In Section \ref{sect:EquiArith}, we use the theorem from the previous section to derive an effective bound for equidistribution along multivariate arithmetic sequences of entries in abelian horospherical flows on the space of lattices.  In this result, we allow the arithmetic sequences in different coordinates to have different spacing, although we will not need it for Section \ref{sect:Sieving}. The techniques used in this section are heavily inspired by Section 3 of \cite{VenkSparse} and also make use of the spectral gap, as well as some Fourier analysis and other analytic techniques.

In Section \ref{sect:Sieving}, we use the bound along arithmetic sequences as well as a combinatorial sieve theorem to obtain an upper and lower bound on averages over almost-prime entries in abelian horospherical flows.  We start with the case of $\Gamma$ cocompact, for which we obtain a result that implies Theorem \ref{thm:CmptThm} above.
We then move to the case $\Gamma=\slnz$,  where we prove a similar result for almost-primes in the orbits of points satisfying a strongly polynomially Diophantine condition,  giving us Theorem \ref{thm:NoncmptThm}.  In order to apply sieving in both cases, we introduce a  particular gcd-sum function (counting the number of integer points in a cube in $\R^d$ of side length $K\in\N$ such that $K$ divides the product of the entries) and verify that the relevant errors can be controlled. 

Finally, in Section \ref{sect:concl}, we make some closing remarks and indicate possible extensions and areas for future research.


\section{Notation and Preliminaries}\label{sect:prelims}

\subsection{Some Basic Notation}

Let $G=\slnr$ for $n\geq2$.  Throughout most of this document, $\Gamma$ will denote  $\slnz$, but we will also discuss the case where $\Gamma\leq G$ is a cocompact lattice. We are interested in the right actions of certain subgroups of $G$ on the right coset space $X=\Gamma\backslash G$.

Although it is not a group, $X$ inherits a finite ``Haar" measure $m_X$ from the (bi-invariant) 
Haar measure $m_G$ on $G$.  In this document, we will always take $m_X$ and $m_G$ to be normalized so that $m_X$ is a probability measure and so that the measure of a small set in $G$ equals the measure of its projection in $X$.
We will use $|\cdot|$ to denote the standard Lebesgue measure on $\R^d$ and $d{\bf t}$ to denote the differential with respect to Lebesgue measure for ${\bf t}\in\R^d$.

We will use gothic letters to represent the Lie algebra of a Lie group (e.g. $\mathfrak{g}$ is the Lie algebra of $G$).  Fix an inner product on $\mathfrak{g}$.  
This extends to a Riemannian metric on $G$ via left translation, which defines a left-invariant metric $d_G$ and a left-invariant volume form, which (by uniqueness) coincides with the Haar measure on $G$ up to scaling.  
This then induces a metric $d_X$ on $X$ of the form
\begin{align*}
d_X(\Gamma g_1, \Gamma g_2) = \inf_{\gamma_1,\gamma_2\in\Gamma} d_G(\gamma_1g_1,\gamma_2g_2) = \inf_{\gamma\in \Gamma} d_G(g_1,\gamma g_2).
\end{align*}
The same construction can be used to define a left-invariant metric $d_H$ for any subgroup $H\leq G$ by restricting the inner product to $\mathfrak{h}\subseteq\mathfrak{g}$.  Note, however, that in general $d_H \neq d_G\vert_H$.  Instead, we have that $d_G(h_1,h_2)\leq d_H(h_1,h_2)$ for $h_1,h_2\in H$, since the infemum used to define the distance $d_G$ is taken over a larger set than in $d_H$.  We will use the notation $B^H_r(h)$ to denote a ball of radius $r$ with respect to the metric $d_H$ around a point $h\in H$ (this is to distinguish these balls from the sets $B_T$ that we will define in Section \ref{sect:horo}).  Also observe that every point has a neighborhood in which the left-invariant metric is Lipschitz equivalent to the metric derived from any matrix norm on ${\rm  Mat}_{n\times n}(\R)$ (see Lemma 9.12 in \cite{EW} for details).

Define the adjoint representation of $g\in G$ as the map ${\rm  Ad}_g: \mathfrak{g}\to\mathfrak{g}$ given by $Y\mapsto gYg^{-1}$ for $Y\in\mathfrak{g}$.

In considering equidistribution questions, our space of test functions will be $C^\infty_c(X)$, the set of smooth, compactly supported (real- or complex-valued) functions on $X$.  Define the action of $G$ on this space by $[g\cdot f](x) = f(xg^{-1})$ for $g\in G$ and $f\in C^\infty_c(X)$.

Finally, we will use the notation $a\ll b$ to indicate that $a$ is less than a fixed constant times $b$ and $a\asymp b$ to indicate that $a\ll b$ and $b\ll a$.  In general, the implied constants may depend on $n$ and on the data of the dynamical system (more specifically, on $d$, the dimension of the horospherical subgroup).
Any additional dependence of the constants will be indicated by a subscript (e.g. $\ll_f$ indicates that the implicit constant may depend on $n$, $d$, and $f$).  In principle, the constants may also depend on the lattice $\Gamma$, although since we are primarily considering $\Gamma=\slnz$, we will not indicate this dependence with a subscript when $\Gamma$ is understood to be fixed in this way.  We will also use the standard notation $\mathcal{O}(f(x))$ to indicate a function whose absolute value is bounded by a constant times $|f(x)|$ as $x\to\infty$, where as before the constant may depend on $n$ and $d$, and any additional dependence will be indicated with a subscript.

\subsection{Horospherical Subgroups}\label{sect:horo}

A subgroup $U$ of $G$ is (expanding) horospherical with respect to an element $g\in G$ if $U=\{u\in G \hspace{3pt}\vert\hspace{3pt} g^{-j}ug^j \to e \text{ as } j\to \infty\}$, where $e$ is the identity.  In other words, elements of $U$ are contracted under conjugation by $g^{-1}$ and expanded under conjugation by $g$. 

Define the one-parameter subgroup $\{a_t\}_{t\in\R}\in G$ by
\begin{align}
a_t &= \exp(t\hspace{2pt}{\rm  diag}(\underbrace{\lambda_1,\cdots, \lambda_1}_{m_1}, \underbrace{\lambda_2,\cdots, \lambda_2}_{m_2}, \cdots, \underbrace{\lambda_N,\cdots, \lambda_N}_{m_N}))\label{eq:a_t}
\end{align}
where $\lambda_1\geq\lambda_2\geq \cdots \geq \lambda_N$.  The requirement that $a_t\in \slnr$ for all $t\in\R$ means that $m_1+\cdots+m_N = n$ and $m_1\lambda_1+\cdots+m_N\lambda_N = 0$.

Let $U$ denote the block-upper-triangular unipotent subgroup given by 

\begin{align}
U=
\left\{\left(\begin{matrix}
\begin{matrix}
I_{m_1}&\vline&&\\
\hline
&\vline&I_{m_2}\\
\end{matrix}&&\bigast\\
&\ddots&\\
\bigzero&&\begin{matrix}
I_{m_{N-1}}&\vline&\\
\hline
&\vline&I_{m_N}
\end{matrix}
\end{matrix} 
\right)\right\}\label{eq:Uhoro}
\end{align}
where $I_m$ is the $m\times m$ identity matrix. Notice that $U$ is the horospherical subgroup corresponding to $a_t$ for $t>0$.  Similarly, define the contracting subgroup $U^{-}$ by
\begin{align*}
U^{-}=
\left\{\left(\begin{matrix}
\begin{matrix}
I_{m_1}&\vline&&\\
\hline
&\vline&I_{m_2}\\
\end{matrix}&&\bigzero\\
&\ddots&\\
\bigast&&\begin{matrix}
I_{m_{N-1}}&\vline&\\
\hline
&\vline&I_{m_N}
\end{matrix}
\end{matrix} 
\right)\right\}
\end{align*}
which is horospherical with respect to $a_{t}$ for $t<0$, and define $U^0$ to be the centralizer of $a_t$ ($t\neq 0$), given by
\begin{align*}
U^{0}=
\left\{\left(\begin{matrix}
\begin{matrix}
B_{1}&\vline&&\\
\hline
&\vline&B_{2}\\
\end{matrix}&&\bigzero\\
&\ddots&\\
\bigzero&&\begin{matrix}
B_{m_{N-1}}&\vline&\\
\hline
&\vline&B_{m_N}
\end{matrix}
\end{matrix} 
\right)\hspace{5pt}\vline \hspace{5pt}
\begin{aligned}
&B_i \in {\rm  GL}_{m_i}(\R)\\
&\det B_{1}\cdots \det B_{N} = 1
\end{aligned}
\right\}.
\end{align*}
Let $d_0 = \sum_{i=1}^N m_i^2$ and observe that $d := \dim U = \dim U^{-} = \frac{1}{2}\left(n^2-d_0\right)$ and $\dim U^0 = d_0 -1$.  All horospherical subgroups of $G=\slnr$ are conjugate to a subgroup of the form given in (\ref{eq:Uhoro}), so we restrict our attention to $U$ of this form.

Observe that $U$ is diffeomorphic to $\R^d$ through any identification ${\bf t}\mapsto u({\bf t})$ of the coordinates of $\R^d$ with the matrix entries in the upper-right corner of (\ref{eq:Uhoro}).\footnote{One could also use the more standard map $u({\bf t}) = \exp(\iota({\bf t}))$, where $\iota:\R^d \mapsto \mathfrak{u}$ is any identification of $\R^d$ with the Lie algebra $\mathfrak{u}$ of $U$.  We have chosen to use the former embedding for ease of notation and because we will later restrict our attention to abelian horosphericals, for which the two maps coincide (up to scaling and permutations of the coordinates).  However, whichever map is used does not substantively change the results presented here.}
Note, however, that $U$ and $\R^d$ are only isomorphic as groups in the case that $U$ is abelian, which occurs when $a_t$ has precisely two eigenvalues.

The bi-invariant Haar measure $m_U$ on $U$ is the pushforward of Lebesgue measure on $\R^d$ under this identification, and we may normalize it so that $u([0,1]^d)$ has unit measure.  Define an expanding family of balls in $U$ by $B_T = a_{\log T} u([0,1]^d) a_{-\log T}$ for $T\in \R$.  One may verify that the preimage of $B_T$ in $\R^d$ is given by a box where $x_k\in[0,T^{\lambda_i-\lambda_j}]$ for $i>j$ if the coordinate $x_k$ is mapped to the $(i,j)$-block of (\ref{eq:Uhoro}) under our identification.  Hence, $m_U(B_T) = T^p$, where $p = \sum_{i>j} m_i m_j (\lambda_i -\lambda_j)$.

\subsection{Measure Decomposition}\label{sect:decomp}

The product map $U\times U^0 \times U^- \to G$ given by $(u, u^0, u^-)\mapsto uu^0u^-$ is a biregular map onto a Zariski open dense subset of $G$ 
(see Proposition 2.7 in \cite{MT}).  In particular, if we let $H = U^0 U^-$, this means that $m_G(G\setminus UH) = 0$ and that the product map $(u, h) \mapsto uh$ is open and continuous.  Additionally, it is not difficult to see that $U\cap H = \{e\}$. 
Then by virtue of the fact that $G$ is unimodular,
we have that $m_G$ restricted to $UH$ is proportional to the pushforward of $m_U \times m^r_H$ by the product map, where $m_H^r$ is the right Haar measure on $H$ (see, e.g., Lemma 11.31 in \cite{EW} or Theorem 8.32 in \cite{Knapp}).  Note that we could equivalently use the left Haar measure on $H$ and multiply by the modular function $\triangle_H$, but for convenience of notation we will use the right Haar measure.

\subsection{Sobolev Norms}\label{sect:Sob}

Fix a basis $\mathcal{B}$ for the Lie algebra $\mathfrak{g}$ of $G$.
Define the (right) differentiation action of $\mathfrak{g}$ on $C^\infty_c(X)$ by $Yf(x) = \frac{d}{dt} f(x \exp(tY))\vert_{t=0}$ for $Y\in \mathcal{B}$ and $f\in C^\infty_c(X)$.  Higher order derivatives of $f$ can then be expressed as monomials in the basis $\mathcal{B}$.

For $p \in [1, \infty]$ and $\ell\in\N$, the $(p,\ell)$-Sobolev norm of $f\in C^\infty_c(X)$ simultaneously controls the $L^p$-norm of all derivatives of $f$ up to order $\ell$.  More precisely, let
\begin{align*}
\sob{p}{\ell}(f) = \sum_{\deg(\mathcal{D})\leq\ell} \nn{\mathcal{D}f}_{L^p(X)}
\end{align*}
where $\mathcal{D}$ ranges over all monomials in $\mathcal{B}$ of degree $\leq \ell$.  Observe that the Sobolev norm can be defined similarly for $C^\infty_c(G)$ and $C^\infty_c(H)$ where $H\leq G$, given a choice of basis for $\mathfrak{h}\subseteq\mathfrak{g}$.\footnote{The choice of the basis $\mathcal{B}$ is unimportant in the sense that choosing a different basis will lead to an equivalent norm.  Likewise, we could use any norm on the components $\nn{\mathcal{D}f}_{L^p(X)}$ (here we have used the $l^1$-norm), but as all such norms are equivalent, the choice is unimportant.}

We will only require the $(2, \ell)$- and $(\infty, \ell)$-Sobolev norms.  When $p=2$, we will drop the notation, letting $\mathcal{S}_\ell(f)=\sob{2}{\ell}(f)$. When needed, we will use a superscript $\mathcal{S}^X$ to indicate a Sobolev norm for functions defined on $X$.

Some useful properties of these norms are as follows (see \cite{VenkSparse} or \cite{KMBddOrbits}):
\begin{enumerate}[(i)]
\item For $X$ a probability space, $f\in C^\infty_c(X)$,
$p\in[1,\infty]$, and $k\leq \ell$, $\sob{p}{k}(f)\leq \sob{\infty}{\ell}(f)$.\label{Sob1}
\item For $f_1, f_2\in C^\infty_c(X)$, $\sob{\infty}{\ell}(f_1f_2)\ll_\ell \sob{\infty}{\ell}(f_1)\sob{\infty}{\ell}(f_2)$.\label{Sob2}
\item For $f\in C^\infty_c(X)$ and $g\in G$, $\sob{\infty}{\ell}(g\cdot f) \ll_\ell \nn{{\rm  Ad}_{g^{-1}}}^\ell \sob{\infty}{\ell}(f)$, where $\nn{\cdot}$ is the operator norm on linear functions $\mathfrak{g}\to\mathfrak{g}$.\label{Sob3}
\item Let $L\subset G$ be compact.  
For $f\in C^\infty_c(X)$, $x\in X$,
\[
|f(xg)-f(x)|\ll_L \sob{\infty}{1}(f)d_G(g,e)
\]
for all $g\in L$.\label{Sob4}
\item Let $X$ and $Y$ be Riemannian manifolds. For $f_1\in C^\infty_c(X)$ and $f_2\in C^\infty_c(Y)$, 
\[\sobt{\ell}^{X\times Y}(f_1 \cdot f_2) \ll_{X,Y} \sobt{\ell}^X(f_1)\sobt{\ell}^{Y}(f_2).
\]\label{Sob5}
\end{enumerate}

\subsection{Approximation to the Identity}\label{sect:approxid}

At times we will want to use smooth bump functions with small support as approximations to the identity, but we will need to know that the Sobolev norm of such functions can be controlled.  For this we have the following lemma, which can be found in \cite{KMBddOrbits}.

\begin{lem}[\hspace{1sp}\cite{KMBddOrbits}, Lemma 2.4.7(b)]\label{lem:approxId}
Let $Y$ be a Riemannian manifold of dimension $k$.  Then for any $0<r<1$ and $y\in Y$, there exists a function $\theta\in C^\infty_c(Y)$ such that:
\begin{enumerate}[(i)]
\item $\theta\geq0$
\item $\text{supp }(\theta)\subseteq B^Y_r(y)$
\item $\int_Y \theta = 1$
\item $\sobt{\ell}^Y(\theta) \ll_{Y,y} r^{-(\ell+k/2)}$.\label{approxidSob}
\end{enumerate}
\end{lem}

\subsection{The Space of Unimodular Lattices}

For $\Gamma = \slnz$, $X$ is noncompact and can be understood as the space of unimodular lattices (that is, lattices of covolume 1) in $\R^n$ under the identification $\Gamma g \leftrightarrow \Z^n g$.

For $0<\epsilon\leq 1$, define $L_\epsilon$ to be the set of lattices in $X = \slnz\backslash\slnr$ with no nonzero vectors shorter than $\epsilon$.  That is, let
\begin{align*}
L_\epsilon = \left\{ \Gamma g \in X \hspace{2pt}\vert\hspace{2pt} \nn{vg} \geq \epsilon\hspace{2pt} \text{ for all } v\in \Z^n \setminus\{0\} \right\}
\end{align*}
where the norm above can be taken to be any norm on $\R^n$, but for convenience we will use the max norm.  By Mahler's Compactness Criterion, $L_\epsilon$ is a compact set (for details and a proof, see \cite{Rag} Corollary 10.9, \cite{BekkaMayer} Theorem 5.3.2, or \cite{EW} Theorem 11.33).

\subsection{Radius of Injection}

Given small $\epsilon>0$, we want to find a radius $r>0$ (depending on $\epsilon$) such that projection at $x$, given by
\begin{align*}
\pi_x: B_r^G(e) &\to B_r^X(x)\\
g &\mapsto xg
\end{align*}
is injective for all $x\in L_\epsilon$ (in fact, it is not difficult to see from the definition of the metric on $X$ that this will be an isometry).  For this, we have the following lemma, which is proved in a much more general setting in \cite{HeeBenoist} (see the proof of Lemma 11.2).  A proof of the lemma as it is stated here can be found in Appendix \ref{app:RadiusofInjection}.

\begin{lem}\label{lem:radiusofinjection}
There exist constants $c_1, c_2>0$ (depending only on $n$) such that for any $0<\epsilon<c_1$, the projection map $\pi_x: B_r^G(e) \to B_r^X(x)$ 
is injective for all $x\in L_\epsilon$, where $r=c_2\epsilon^n$.
\end{lem}

\subsection{Quantitative Nondivergence}\label{sect:nondiv}

Let $\{e_1, \cdots, e_n\}$ be the standard basis on $\R^n$.  Let $e_I = e_{i_1}\wedge\cdots\wedge e_{i_j}$ for a multi-index $I=(i_1, \cdots, i_j)$, where $1\leq i_1<\cdots<i_j\leq n$.   Then $\{e_I \}$ is a basis for $\Lambda^j(\R^n)$, the $j$th exterior power of $\R^n$.  Define the norm of $w=\sum_{I} w_I e_I \in\Lambda^j(\R^n)$ to be $\nn{w}=\max_I |w_I|$.  Denote by $\Lambda^j(\Z^n)$ the discrete subset of $\Lambda^j(\R^n)$ composed of linear combinations of basis vectors with integer coefficients. Notice that $g\in{\rm  GL}_n(\R)$ acts on $\Lambda^j(\R^n)$ on the right by
\begin{align*}
(e_{i_1}\wedge\cdots\wedge e_{i_j})g = (e_{i_1}g)\wedge\cdots\wedge (e_{i_j}g)
\end{align*}
where the action extends to all of $\Lambda^j(\R^n)$ via linearity.

The following theorem quantitatively describes how often certain polynomial maps from $\R^d$ to $X$ land inside a compact set $L_\epsilon$.  This is a special case of Theorem 5.2 in \cite{KMNondiv}, which itself extends results of \cite{DaniNondiv} and \cite{MargulisNondiv}.  The original theorem is stated for much more general $(C,\alpha)$-good functions, but we will only need the version below, which uses the observation in Lemma 3.2 of \cite{BKMNondiv}
that polynomials in $\R[x_1,\cdots,x_d]$ of degree $\leq k$ are $(C_{d,k}, 1/dk)$-good on $\R$.

\begin{thm}[\hspace{1sp}\cite{KMNondiv}, Theorem 5.2]\label{thm:nondiv1}
Let $d, n, k\in\N$ and $0<\rho\leq1/n$.  Let $B\subset\R^d$ be a ball and suppose $\dst{\xi:B\to {\rm GL}_n(\R)}$ satisfies:
\begin{enumerate}[(i)]
\item $\nn{w\xi({\bf t})}$ is a polynomial in the coordinates of ${\bf t}$ of degree $\leq k$, and\label{nondiv1}
\item $\dst{\sup_{{\bf t}\in B} \nn{w\xi({\bf t})} \geq \rho}$\label{nondiv2}
\end{enumerate}
for all primitive $w\in \Lambda^j(\Z^n)\setminus\{0\}$ and $j\in\{1,\cdots, n\}$.  Then for any $0<\epsilon\leq\rho$, 
\begin{align*}
|\{{\bf t}\in B \hspace{2pt} \vert \hspace{2pt} \Gamma\xi({\bf t}) \notin L_\epsilon \} | \ll_{d,k} \left(\epsilon/\rho\right)^{1/dk} |B|.
\end{align*}
\end{thm}

From this theorem we may derive the following corollary, which we will use in the proof of Theorem \ref{thm:equidist} to say that the orbit of a point satisfying a certain Diophantine condition spends a relatively large proportion of time in $L_\epsilon$ when pushed by the flow $a_t$.

\begin{cor}\label{cor:nondiv}
Let $T,R>1$ and $x_0=\Gamma g_0 \in X$.  Then suppose $R_0>0$ is such that
\begin{align*}
\sup_{{\bf t}\in [0,1]^d} \nn{wg_0a_{\log T} u({\bf t}) a_{-\log T}}\geq R_0
\end{align*}
for all primitive $w\in \Lambda^j(\Z^d)\backslash\{0\}$ and $j\in\{1,\cdots,n-1\}$ and define $\rho = \min(1/n, R_0/R^q)$.  Then for any $0<\epsilon<\rho$,
\begin{align*}
\left| \{ {\bf t}\in [0,1]^d \hspace{2pt}\vert\hspace{2pt} x_0 a_{\log T} u({\bf t}) a_{-\log T}a_{\log R} \notin L_\epsilon \} \right|\ll (\epsilon/\rho)^{1/d(n-1)}.
\end{align*}
\end{cor}

\begin{proof}
Let $\xi({\bf t}) = g_0 a_{\log T} u({\bf t}) a_{-\log T}a_{\log R}$.  We want to demonstrate that conditions (\ref{nondiv1}) and (\ref{nondiv2}) hold in Theorem \ref{thm:nondiv1} for $k=n-1$, $\rho = \min(1/n, R_0/R^q)$, and $B=[0,1]^d$.

Recall that our identification $u({\bf t})$ places one coordinate of ${\bf t}$ in each matrix entry in the upper-right corner of (\ref{eq:Uhoro}).  Then since multiplication by $a_t$ on either the left or the right only changes matrix entries by scaling, each entry in the upper-right corner of $a_{\log T} u({\bf t}) a_{-\log T}a_{\log R}$ only depends linearly on a single coordinate of ${\bf t}$.  This means that for any matrix $g_0$, all entries of  $\xi({\bf t}) = g_0 a_{\log T} u({\bf t}) a_{-\log T}a_{\log R}$ will be affine.
Hence, when we take wedge products of the form
\begin{align*}
(e_{i_1}\wedge\cdots\wedge e_{i_j})\xi({\bf t}) = (e_{i_1}\xi({\bf t}))\wedge\cdots\wedge (e_{i_j}\xi({\bf t}))
\end{align*}
the coefficients will be polynomials of degree $\leq j$.  Furthermore, since $\xi({\bf t})\in\slnr$ for all ${\bf t}$, $(e_1\wedge\cdots\wedge e_n)\xi({\bf t})=e_1\wedge\cdots\wedge e_n$ is independent of ${\bf t}$, and the top exterior power can be ignored.  Then from the definition of the norm on $w\in\Lambda^j(\R^n)$, we have that $\nn{w\xi({\bf t})}$ is a polynomial of degree $\leq n-1$ for all $j\in\{1,\cdots, n\}$, so (\ref{nondiv1}) is satisfied with $k=n-1$.

Moreover, notice that $e_k a_{\log R}=R^{\lambda_i}e_k$ if $\lambda_i$ is $k$\textsuperscript{th} eigenvalue in the definition of $a_t$ in (\ref{eq:a_t}).  Then the right action of $a_{\log R}$ scales $e_{i_1}\wedge\cdots\wedge e_{i_j}\in\Lambda^j(\R^d)$ by the product of all such corresponding factors. 
Since $R>1$, the most $a_{\log R}$ can therefore contract any basis element is by the product of all scaling factors corresponding to negative eigenvalues of (\ref{eq:a_t}), that is, by $R^{-q}$, where $q = \sum_{\lambda_i <0} -m_i\lambda_i$.  It then follows from the definition of the norm that 
\begin{align}
\nn{wa_{\log R}}\geq R^{-q}\nn{w}\label{eq:R-q}
\end{align}
for any $w\in\Lambda^j(\R^d)\setminus\{0\}$ and $j\in\{1,\cdots, n\}$.

Now observe that for $\rho = \min(1/n, R_0/R^q)$, we have $0<\rho\leq 1/n$ and also
\begin{align*}
\sup_{{\bf t}\in[0,1]^d}\nn{w\xi({\bf t})} &= \sup_{{\bf t}\in[0,1]^d}\nn{wg_0 a_{\log T} u({\bf t}) a_{-\log T}a_{\log R}}\\
&\geq R^{-q}\sup_{{\bf t}\in[0,1]^d}\nn{wg_0 a_{\log T} u({\bf t}) a_{-\log T}}\\
&\geq R_0/R^q\\
&\geq \rho
\end{align*}
for $j\in\{1,\cdots,n-1\}$ and primitive $w\in\Lambda^j(\Z^n)\setminus\{0\}$.  Thus condition (\ref{nondiv2}) is satified, since as before, $\xi({\bf t}) \in \slnr$ implies the condition is trivially satisfied for the top exterior power.

Hence, by Theorem \ref{thm:nondiv1}, we have
\begin{align*}
\left| \{ {\bf t}\in [0,1]^d \hspace{2pt}\vert\hspace{2pt} \Gamma \xi({\bf t}) \notin L_\epsilon \} \right| \ll (\epsilon/\rho)^{1/d(n-1)}.
\end{align*}
\end{proof}

\subsection{Decay of Matrix Coefficients}\label{sect:expmixing}

In order to obtain effective rates of equidistribution in Sections \ref{sect:Equidistribution} and \ref{sect:EquiArith}, we will need to use results on the effective decay of matrix coefficients.

Estimates of this type have a long and rich history, including Selbrerg's celebrated $3/16$ theorem for congruence quotients of $\slz$, Kazhdan's property (T), and works of Cowling, Moore, Howe, and Oh. Far reaching extensions of Selberg's work are also in place thanks to works of Jacques-Langlands, Burger-Sarnak, and Clozel. Our formulation here is taken from \cite{KMBddOrbits} (see \cite{KMBddOrbits}, \cite{GMO}, and \cite{EMMV} for a more comprehensive history and discussion).

\begin{thm}[\hspace{1sp}\cite{KMBddOrbits}, Corollary 2.4.4]
Let $G=\slnr$ and $X=\Gamma/G$ for a lattice $\Gamma$.  There exists a constant $0<\beta<1$ such that for $f_1, f_2 \in \CX$ and $g\in G$, 
\begin{align*}
\left|\left< g\cdot f_1, f_2\right>_{L^2(X)} - \int_X f_1 dm_X \int_X \overline f_2 dm_X \right|\ll e^{-\beta d_G(e,g)} \sobt{\ell}(f_1)\sobt{\ell}(f_2)
\end{align*}
where $\ell$ is the dimension of maximal compact subgroup of $G$.  When $n\geq 3$, the constant $\beta$ is independent of the lattice $\Gamma$, and when $n=2$ it is independent of the lattice if $\Gamma$ is a congruence lattice.
\end{thm}

For our specific applications, we have the following immediate corollaries.

\begin{cor}\label{cor:expmixing}
Let the setting be as above.
\begin{enumerate}[(i)]
\item For $f_1, f_2\in\CX$ and $t\geq 0$, we have
\[
\left|\int_X f_1(xa_t)f_2(x)dm_X(x)-\int_X f_1 dm_X\int_X \overline f_2 dm_X\right|\ll e^{-\beta t}\sobt{\ell}(f_1)\sobt{\ell}(f_2)
\]\label{expmixing}
\item For $f\in\CX$ and ${\bf t}\in \R^d$, 
\[
\left|\left<u({\bf t})f , f\right>_{L^2(X)}-\left|\int_X f dm_X\right|^2\right|\ll\max(1,|{\bf t}|)^{-\beta}\sobt{\ell}(f)^2
\]\label{umtxcoeff}
\end{enumerate}
\end{cor}

\subsection{Combinatorial Sieve}\label{sect:Sieve}

In order to understand the distribution of almost-prime times in horospherical orbits we will make use of the following combinatorial sieve theorem (see \cite{CombSieve}, or \cite{SarnackNevo} for a form more similar to that stated here).

\begin{thm}[\hspace{1sp}\cite{CombSieve}, Theorem 7.4]\label{thm:CombSieve}
Let $A = \{a_n\}$ be a sequence of nonnegative numbers and let $\dst{P = P(z) = \prod_{p<z} p}$ be the product of primes less than $z$.  Let $\dst{S(A, P) = \sum_{(n,P)=1} a_n}$ and $\dst{S_K(A,P) = \sum_{ n \equiv 0 \mod K} a_n}$. 
Then suppose
\begin{enumerate}[(i)]
\item \label{axiom1} There exists a multiplicative function $g(K)$ on $K$ squarefree such that
\begin{align*}
S_K(A,P) = g(K)\mathcal{X} + r_K(A)
\end{align*}
and for some $c_1 > 0$, we have $0\leq g(p) < 1 - \frac{1}{c_1}$ for all primes $p$. 
\item \label{axiom2} $A$ has level distribution $D(\mathcal{X})$, i.e. there is $\epsilon > 0$ such that
\begin{align*}
\sum_{K<D} |r_K(A)| \ll_\epsilon \mathcal{X}^{1-\epsilon}.
\end{align*}
\item \label{axiom3} $A$ has sieve dimension $r$, i.e. there exists a constant $c_2>0$ such that for all $2 \leq w \leq z$, we have
\begin{align*}
-c_2 \leq \sum_{w \leq p \leq z} g(p)\log p - r \log \frac{z}{w} \leq c_2.
\end{align*}
\end{enumerate}
Then for $s>9r$, $z=D^{1/s}$, and $\mathcal{X}$ large enough, we have
\begin{align*}
S(A,P) \asymp \frac{\mathcal{X}}{(\log\mathcal{X})^r}
\end{align*}
where the implicit constants depend on the constants in (\ref{axiom1}), (\ref{axiom2}), and (\ref{axiom3}).
\end{thm}


\section{Effective Equidistribution of Horospherical Flows}\label{sect:Equidistribution}

Our main objective in this section is to prove the following effective equidistribution theorem for horospherical flows on $X=\slnz\backslash\slnr$.

\begin{thm}\label{thm:equidist}
Let $u$ be a horospherical flow on $X$, and let $x_0=\Gamma g_0\in X$.  Then there exist constants $\gamma, C>0$ (depending only on $n$ and $d$)
such that for $f\in C^\infty_c(X)$ and $T>R>C$, either
\begin{align}
\left| \frac{1}{m_U(B_T)} \int_{B_T} f(x_0 u) dm_U(u) - \int_X f dm_X \right| \ll R^{-\gamma} \sob{\infty}{\ell}(f) \tag{\ref{thm:equidist}.a}\label{eq:endresult}
\end{align}
or
\begin{align}
\exists j\in\{1,\cdots, n-1\} \text{ and primitive } w\in\Lambda^j(\Z^n)\setminus\{0\}\text{ s.t. } \nn{w g_0 u}< R^q\hspace{5pt}\forall {u}\in B_T.\tag{\ref{thm:equidist}.b}\label{eq:dionegation}
\end{align}
where $q = \sum_{\lambda_i <0} -m_i\lambda_i$ and $\ell=n(n-1)/2$ is the dimension of maximal compact subgroup of $G$.
\end{thm}

Intuitively, this theorem says that either the $U$-orbit of $x_0$ equidistributes in $X$ with a fast rate, or $x_0$ is close to a proper subset of $X$ that is fixed by the action of $U$, where our notions of ``fast" and ``close" are quantitatively related.

\begin{rmk}
The condition that $w$ in (\ref{eq:dionegation}) be primitive is conceptually useful but technically  unnecessary, in that if there exists any $w\in\Lambda^j(\Z^n)\setminus\{0\}$ satisfying (\ref{eq:dionegation}), then there will also exist a primitive vector that does so.
\end{rmk}

\begin{rmk}
The ``either/or" in the theorem statement is not meant to imply an exclusive or.  In fact, the theorem can be restated in the following form: For $x_0$, $\gamma$, $C$, $\ell$, $f$, $T$, and $R$ as above, not (\ref{eq:dionegation}) implies (\ref{eq:endresult}).

This leads us to define the following Diophantine basepoint condition for $x_0 = \Gamma g_0 \in X$:
\begin{align}
\forall j\in\{1,\cdots,n-1\} \text{ and } w\in\Lambda^j(\Z^n)\setminus\{0\}, \hspace{4pt}\exists\hspace{2pt} u\in B_T \text{ s.t. } \nn{w g_0u}\geq R^q.\tag{\ref{thm:equidist}.c}\label{eq:dio1}
\end{align}
Then Theorem \ref{thm:equidist} says that (\ref{eq:dio1}) implies (\ref{eq:endresult}), and this is in fact how we will structure the proof.
\end{rmk}

\begin{rmk}
Although the theorem is stated for balls of the form $B_T = a_{\log T} u([0,1]^d) a_{-\log T}$, it holds equally well for symmetric balls of the form $B_T = a_{\log T} u([-1,1]^d) a_{-\log T}$.
\end{rmk}

\begin{proof}
Let $x_0 = \Gamma g_0\in X$ satisfy the basepoint condition in (\ref{eq:dio1}) for some $T>R$.  Then consider $f\in C_c^\infty(X)$ and write, via a change of variables,
\begin{align}
I_0 &:= \frac{1}{m_U(B_T)}\int_{B_T} f(x_0u) dm_U(u)\nonumber\\
&= \frac{1}{m_U(B_{T/R})} \int_{B_{T/R}} f(x_0a_{\log R} u a_{-\log R})dm_U(u)\nonumber\\
&= \frac{1}{m_U(B_{T/R})} \int_{U} {\bf1}_{B_{T/R}}(u) f(x_0a_{\log R} u a_{-\log R})dm_U(u).\label{eq:I0}
\end{align}
We want to show that this quantity is close to $\int f dm_X$, and from (\ref{eq:I0}) it almost looks as if we could apply the exponential mixing result of Corollary \ref{cor:expmixing} (\ref{expmixing}) to achieve this, however there are several significant barriers to doing so.  Most obviously, the integral in (\ref{eq:I0}) is over $U$ instead of $X$.  Furthermore, the ``basepoint" $x_0a_{\log R} u$ varies with $u$, and will eventually spend time outside of any fixed compact subset of $X$ for $u$ coming from a large enough ball. Finally, the function ${\bf1}_{B_{T/R}}$ is not smooth.

We will first address the issue of smoothness by convolving the indicator function with a smooth approximation to the identity ({\bf Step 1}).  We will then apply the ``thickening" argument of Margulis to obtain an intergral over $X$ from our inegral over $U$ ({\bf Step 2}).  Finally, we will deal with the moving basepoint by demonstrating that for most $u\in B_{T/R}$ we have a uniformly good rate of equidistribution and that the size of the set on which this does not occur can be quantitatively controlled ({\bf Step 3}).  This last step is where we will use the nondivergence result of Section \ref{sect:nondiv}. 

\subsection*{Step 1}
Let $r$ be a small, positive number (to be determined) and let $\theta\in C_c^\infty (U)$ be a nonnegative bump function supported on $B^U_r(e)$ satisfying the approximate identity properties of Lemma \ref{lem:approxId}.  Then the convolution $\int_U\theta(u'){\bf1}_{B_{T/R}}(u(u')^{-1})dm_U(u')$ is a smooth function approximating our original indicator function.  If we substitute this function for ${\bf1}_{B_{T/R}}$ in (\ref{eq:I0}) and use the invariance property of the Haar measure, we get the integral
\begin{align}
\Ismth &:= \frac{1}{m_U(B_{T/R})} \int_{U} \int_U\theta(u'){\bf1}_{B_{T/R}}(u(u')^{-1})dm_U(u') f(x_0a_{\log R} u a_{-\log R})dm_U(u)\nonumber\\
&=\frac{1}{m_U(B_{T/R})} \int_{U}\int_U \theta(u') {\bf 1}_{B_{T/R}}(u) f(x_0a_{\log R} uu' a_{-\log R})dm_U(u)dm_U(u').\label{eq:Ismth}
\end{align}

Now observe that since $\int\theta =1$, we may again use the invariance of the Haar measure to rewrite (\ref{eq:I0}) as
\begin{align}
I_0 &= \frac{1}{m_U(B_{T/R})} \int_{U}  {\bf 1}_{B_{T/R}}(uu')f(x_0a_{\log R} uu' a_{-\log R})dm_U(u) \int_U \theta(u') dm_U(u')\nonumber\\
&= \frac{1}{m_U(B_{T/R})} \int_{U}\int_U \theta(u') {\bf 1}_{B_{T/R}}(uu') f(x_0a_{\log R} uu' a_{-\log R})dm_U(u)dm_U(u').\label{eq:shiftedBTR}
\end{align}
From (\ref{eq:Ismth}) and (\ref{eq:shiftedBTR}), we can see that
\begin{align}
\left| I_0  - \Ismth\right| &\leq \frac{1}{m_U(B_{T/R})}\int_U\theta(u') \sob{\infty}{0}(f)\left(\int_U \left|{\bf 1}_{B_{T/R}}(uu')-{\bf1}_{B_{T/R}}(u)\right|dm_U(u)\right)dm_U(u')\nonumber\\
&= \frac{\sob{\infty}{0}(f)}{m_U(B_{T/R})}\int\theta(u') m_U(B_{T/R}\triangle B_{T/R}(u')^{-1}) dm_U(u').\label{eq:symmdiffint}
\end{align}
But notice that since $\text{supp }\theta \subseteq B^U_r(e)$, we know that $u'$ is close to the identity, so $u$ in this region can only shift $B_{T/R}$ by a small amount.  In fact, by pulling the measure back to $\R^n$, one may compute directly that the size of the symmetric difference is bounded by
\begin{align}
m_U(B_{T/R}\triangle B_{T/R}(u')^{-1}) \ll (T/R)^{p-p_0}r\label{eq:symmdiffU}
\end{align}
for any $u'\in B^U_r(e)$, where $p_0 = \min_{i>j} (\lambda_i - \lambda_j)$.
Combining this with (\ref{eq:symmdiffint}) above and again using the fact that $\theta$ integrates to 1, we see that
\begin{align}
|I_0-\Ismth| &\ll \frac{(T/R)^{p-p_0}}{m_U(B_{T/R})}r \sob{\infty}{0}(f)=(R/T)^{p_0} r \sob{\infty}{0}(f) \leq r \sob{\infty}{0}(f)\label{eq:I0-Ismth}
\end{align}
since $m_U(B_{T/R}) = (T/R)^p$ and $T\geq R$.

Now that we know $I_0$ and $\Ismth$ can be made close, we want to know that $\Ismth$ is not too far from $\int f dm_X$.  Using Fubini's Theorem, we can say
\begin{align*}
\Ismth = \frac{1}{m_U(B_{T/R})} \int_{B_{T/R}}\int_U \theta(u') f(x_0a_{\log R} uu' a_{-\log R})dm_U(u')dm_U(u)
\end{align*}
and we may also write
\begin{align*}
\int_X f dm_X = \frac{1}{m_U(B_{T/R})} \int_{B_{T/R}} \left( \int_X f dm_X \right) dm_U(u).
\end{align*}
Hence, 
\begin{align}
&\left| \Ismth - \int_X f dm_X \right|\label{eq:mixinginint}\\
&\leq \frac{1}{m_U(B_{T/R})} \int_{B_{T/R}} \left| \int_U \theta(u') f(x_0a_{\log R} uu' a_{-\log R})dm_U(u') -\int_X f dm_X\right| dm_U(u).\nonumber
\end{align}

\subsection*{Step 2}
Now the expression inside the absolute value looks more similar to that of Corollary \ref{cor:expmixing} (\ref{expmixing}), but we are still integrating over the wrong space.  We want an integral over $X$, and although functions on $X$ integrate locally like their pullback by projection over $G$, the integral with which we are concerned is over the lower-dimensional (``thin") subspace $U$.

Define
\begin{align}
I_U(u) := \int_U \theta(u') f(x_0a_{\log R} uu' a_{-\log R})dm_U(u').\label{eq:Ithin}
\end{align}
to be the integral from inside (\ref{eq:mixinginint}) above.  In order to apply exponential mixing, we will need to ``thicken" this integral over $U$ to an integral over a neighborhood of the orbit in $G$ and then project to $X$. 

Recall from Section \ref{sect:decomp} that $m_G = m_U \times m^r_H$, where $m_H^r$ is the right Haar measure on $H = U^0 U^-$.  Then let $\psi\in C^\infty_c(H)$ be an approximate identity supported on $B^H_r(e)$ as described in Lemma \ref{lem:approxId}.
Since $\int \psi = 1$, we may rewrite (\ref{eq:Ithin}) as
\begin{align}
I_U(u) = \int_H \int_U \theta(u') \psi(h) f(x_0a_{\log R} uu' a_{-\log R}) dm_U(u') dm^r_H(h).
\end{align}
Now define 
\begin{align}
I_X(u) := \int_H \int_U \theta(u') \psi(h) f(x_0a_{\log R} uu'h a_{-\log R}) dm_U(u')dm^r_H(h)\label{eq:Ithick}
\end{align}
which differs from $I_U(u)$ only by the presence of the variable $h$ inside $f$.  To see that $I_U(u)$ and $I_X(u)$ are close, observe that
\begin{align}
\left|I_U(u) - I_X(u)\right| \leq \int_H \int_U \theta(u') \psi(h) \left|f(\tilde x) - f(\tilde x a_{\log R} h a_{-\log R})\right| dm_U(u')dm^r_H(h)\label{eq:thin-thick1}
\end{align}
where $\tilde x = x_0a_{\log R} uu'a_{-\log R}$.  But since $f$ has bounded derivative,
\begin{align}
\left|f(\tilde x) - f(\tilde x a_{\log R} h a_{-\log R})\right| \ll \sob{\infty}{1}(f) d_G(e, a_{\log R} h a_{-\log R})\label{eq:bddderiv}
\end{align}
by Sobolev property (\ref{Sob4}).  Furthermore, since conjugation by $a_t$ is non-expanding on the subgroup $H$ (recall that it fixes $U^0$ and contracts $U^-$), we may see
 that\footnote{There is a slight subtlety here because we used the right Haar measure on $H$, so the corresponding metric $d_H$ is right-invariant, while $d_G$ is left-invariant.  In general, $d_G$ restricted to $H$ will be less than or equal to the corresponding left-invariant metric on $H$.  However, any left-invariant metric is Lipschitz equivalent to any right-invariant metric in a suitable neighborhood of the identity, so the above series of inequalities goes through for $r$ small enough.}
\begin{align}
d_G(e, a_{\log R}h a_{-\log R})\leq d_G(e, h)\ll d_H(e,h)\leq r\label{eq:Hnonexpanding}
\end{align}
for $h\in \text{supp } \psi \subseteq B^H_r(e)$.

Then from (\ref{eq:thin-thick1}), (\ref{eq:bddderiv}), and (\ref{eq:Hnonexpanding}) and the fact that both $\theta$ and $\psi$ integrate to 1, we have 
\begin{align}
\left|I_U(u) - I_X(u)\right| &\ll
\sob{\infty}{1}(f) r.\label{eq:thin-thick}
\end{align}

Now we want to verify that $I_X(u)$ is not far from $\int f dm_X$. By our measure decomposition, we can see (\ref{eq:Ithick}) as an integral over $G$:
\begin{align}
I_X(u) = \int_G \phi(g) f(x_0a_{\log R} ug a_{-\log R}) dm_G(g)\label{eq:intoverG}
\end{align}
where the function $\phi(uh) = \theta(u)\psi(h)$ is defined for all $g\in UH$, hence it is defined almost-everywhere.  In order to apply mixing, we want to further interpret $I_X(u)$ as an integral over $X$.  To do this, let $y=x_0 a_{\log R} u$, keeping in mind that $y$ depends on $u$.  Then define $\phi_y\in C^\infty_c(X)$ by $\phi_y = \phi \circ \pi_y^{-1}$ where $\pi_y: G \to X$ is natural projection at $y$.  Note, however, that $\phi_y$ is only well-defined if $\pi_y$ is injective on $\text{supp }\phi = \text{supp }\theta\text{ supp }\psi \subseteq B^U_r(e)B^H_r(e)$.  In a neighborhood of the identity, $B^U_r(e)B^H_r(e) \subseteq B^G_{cr}(e)$ for a positive constant $c$, since
\begin{align*}
d_G(uh,e)\leq d_G(uh,u)+d_G(u,e) =  d_G(h,e)+d_G(u,e) \ll d_H(h,e) + d_U(u,e) \leq 2r.
\end{align*}
Therefore, if $\pi_y$ is injective on $B^G_{cr}(e)$ for $y=x_0a_{\log R}u$ (an assumption we will reutrn to later) we can say from (\ref{eq:intoverG}) that
\begin{align*}
I_X(u) &= \int_G \phi(g) f(yg a_{-\log R}) dm_G(g)\\
&=\int_G \phi_y(yg) f(yg a_{-\log R}) dm_G(g)\\
&= \int_X \phi_y(x)f(xa_{-\log R}) dm_X(x).
\end{align*}
Since $\int \phi_y dm_X = \int \phi dm_G = \int \theta dm_U \int \psi dm^r_H = 1$, we can now apply the effective mixing result from Section \ref{sect:expmixing} to obtain
\begin{align*}
\left|I_X(u) - \int_X f dm_X \right| &= \left|\int_X \phi_y(x)f(xa_{-\log R}) dm_X(x) - \int_X \phi_y dm_X \int_X f dm_X \right|\\
&\ll  R^{-\beta} \sobt{\ell}(\phi_y)\sobt{\ell}(f).
\end{align*}
Then from property (\ref{Sob5}) in Section \ref{sect:Sob} and our bound on the Sobolev norm of an approximate identity (property (\ref{approxidSob}) in Section \ref{sect:approxid}), we can say 
\begin{align*}
\sobt{\ell}^X(\phi_y) = \sobt{\ell}^G(\phi) \ll \sobt{\ell}^U(\theta)\sobt{\ell}^H(\psi) \ll r^{-(\ell+d/2)}r^{-(\ell+\tilde d/2)} = r^{-2\ell - (n^2-1)/2}
\end{align*}
where $\tilde d = \dim H$.
Thus if $\pi_y$ is injective on $B^G_{cr}(e)$, then
\begin{align}
\left|I_X(u) - \int_X f dm_X \right| \ll  R^{-\beta}  r^{-p_1} \sobt{\ell}(f)\label{eq:thick-int}
\end{align}
where $p_1 = 2\ell+(n^2-1)/2$.

\subsection*{Step 3}
However, as we have noted, $y=x_0a_{\log R}u$ depends on $u$, which varies over $B_{T/R}$ in (\ref{eq:mixinginint}).  While we cannot ensure that $\pi_y$ is injective on $B^G_{cr}(e)$ for all $u\in B_{T/R}$, we can say that the set on which this does not occur has small measure.

Recall from Lemma \ref{lem:radiusofinjection} that $\pi_y: B_r^G(e) \to B_r^X(y)$ is injective for $y\in L_\epsilon$ for $r$ proportional to $\epsilon^n$ and for $\epsilon$ small enough.  Furthermore, observe that condition (\ref{eq:dio1}) is equivalent to the statement that for all $j\in\{1, \cdots, n-1\}$ and primitive $w\in\Lambda^j(\Z^d)\setminus\{0\}$, there exists ${\bf t}\in[0,1]^d$ such that $\nn{wg_0a_{\log T} u({\bf t}) a_{-\log T}}\geq R^q$.  Then by Corollary \ref{cor:nondiv} in Section \ref{sect:nondiv}, we have that 
\begin{align*}
\left| \{ {\bf t}\in [0,1]^d \hspace{2pt}\vert\hspace{2pt} x_0 a_{\log T} u({\bf t}) a_{-\log T}a_{\log R} \notin L_\epsilon \} \right|\ll \epsilon^{1/d(n-1)}.
\end{align*}
since $R_0 = R^{-q}$ implies $\rho = 1/n$.  
From this we find that
\begin{align*}
\left| \{ {\bf t}\in [0,1]^d \hspace{2pt}\vert\hspace{2pt} x_0 a_{\log T} u({\bf t}) a_{-\log T}a_{\log R} \notin L_\epsilon \} \right|
&= \left| \{ {\bf t}\in [0,1]^d \hspace{2pt}\vert\hspace{2pt} x_0 a_{\log R} a_{\log T/R} u({\bf t}) a_{-\log T/R} \notin L_\epsilon \} \right|\\
&= m_U\left(\{ u\in B_1 \hspace{2pt}\vert\hspace{2pt} x_0 a_{\log R} a_{\log T/R} u a_{-\log T/R} \notin L_\epsilon \}\right)\\
&= m_U\left(\{ u\in B_{T/R} \hspace{2pt}\vert\hspace{2pt} x_0 a_{\log R} u \notin L_\epsilon \}\right)/m_U(B_{T/R})
\end{align*}
where the last equality can be verified using a change of variables.  That is, for $x_0$ satisfying condition (\ref{eq:dio1}), we have
\begin{align*}
m_U\left(\{ u\in B_{T/R} \hspace{2pt}\vert\hspace{2pt} x_0 a_{\log R} u \notin L_\epsilon \}\right) \ll \epsilon^{1/d(n-1)} m_U(B_{T/R}).
\end{align*}
In other words, if we let $E:= \{ u\in B_{T/R} \hspace{2pt}\vert\hspace{2pt} x_0 a_{\log R} u \in L_\epsilon \}$, then (\ref{eq:thick-int}) holds for all $u\in E$ and $m_U(B_{T/R}\setminus E) \ll \epsilon^{1/d(n-1)} m_U(B_{T/R})$. Thus, from (\ref{eq:mixinginint}), (\ref{eq:thin-thick}), and (\ref{eq:thick-int}), we find
\begin{align*}
\left| \Ismth - \int_X f dm_X \right| \leq& \frac{1}{m_U(B_{T/R})} \int_{B_{T/R}} \left| I_U(u) -\int_X f dm_X\right| dm_U(u)\\
\leq& \frac{1}{m_U(B_{T/R})} \int_{B_{T/R}} \left| I_U(u) -I_X(u) \right| dm_U(u)\\
&\hspace{5pt}+ \frac{1}{m_U(B_{T/R})} \int_{B_{T/R}}\left|I_X(u) - \int_X f dm_X\right| dm_U(u)\\
\ll& \sob{\infty}{1}(f)r +\frac{1}{m_U(B_{T/R})} \int_E \left|I_X(u) - \int_X f dm_X\right| dm_U(u)\\
&\hspace{5pt}+ \frac{1}{m_U(B_{T/R})} \int_{B_{T/R}\setminus E} \left|I_X(u) - \int_X f dm_X\right| dm_U(u)\\
\ll& \sob{\infty}{1}(f)r + \frac{m_U(E)}{m_U(B_{T/R})} R^{-\beta}  r^{-p_1} \sobt{\ell}(f) + \frac{m_U(B_{T/R}\setminus E)}{m_U(B_{T/R})} \sob{\infty}{0}(f)\\
\ll& \sob{\infty}{1}(f)r +R^{-\beta}  r^{-p_1} \sobt{\ell}(f) + \epsilon^{1/d(n-1)} \sob{\infty}{0}(f).
\end{align*}
Finally, from this and (\ref{eq:I0-Ismth}), we have
\begin{align*}
\left| I_0 - \int_X f dm_X \right| &\leq \left| I_0 - \Ismth\right| + \left| \Ismth - \int_X f dm_X \right|\\
&\ll \left(\epsilon^n + R^{-\beta}\epsilon^{-p_1 n} + \epsilon^{1/d(n-1)}\right)\sob{\infty}{\ell}(f)
\end{align*}
where we have used that $r$ is proportional to $\epsilon^n$, as well as Sobolev property (\ref{Sob1}).

Let $p_2:= 1/d(n-1)$. Since $n>p_2$, the $\epsilon^n$ term above decays more quickly than other terms and can be ignored.  To optimize the rate of decay, we set
\begin{align*}
R^{-\beta}\epsilon^{-p_1 n} = \epsilon^{p_2}
\end{align*}
which implies
\begin{align*}
\epsilon = R^{-\beta/(p_1 n + p_2)}.
\end{align*}
Then so long as $R$ is chosen sufficiently large so that $\epsilon$ (and subsequently $r$) are small enough to make Corollary \ref{cor:nondiv} and Lemma \ref{lem:radiusofinjection} true (along with several other statements we made regarding neighborhoods of the identity), then
we have demonstrated (\ref{eq:endresult}) in Theorem \ref{thm:equidist} with the rate
\begin{align*}
\left| I_0 - \int_X f dm_X \right| \ll R^{-\gamma}\sob{\infty}{\ell}(f)
\end{align*}
where $\gamma = \beta p_2/(p_1n+p_2)$.

\end{proof}

\begin{rem}
In the case of $\Gamma$ cocompact, it follows from the above proof that we may remove dependence on the basepoint from our effective equidistribution statement.  That is, for $X=\Gamma\backslash G$, $\Gamma\leq G$ a cocompact lattice, and $U\leq G$ a horospherical subgroup, we have that there exists $\gamma>0$ 
(depending\footnote{ Since dependence on $\Gamma$ only arises from the spectral gap, we can remove dependence for $n\geq 3$ or for $n=2$ when $\Gamma$ is a congruence lattice.}
only on $n$, $d$, and $\Gamma$) such that for $T$ large enough,
\begin{align}
\left|\frac{1}{m_U(B_T)}\int_{B_T} f(x_0 u)dm_U(u) - \int_X f dm_X\right|\ll_\Gamma T^{-\gamma}\sob{\infty}{\ell}(f)\label{eq:equidistcmpt}
\end{align}
for any $f\in C^\infty(X)$ and $x_0\in X$.  This is because we only make use of the basepoint condition in {\bf Step 3}, where we need it to deal with the moving basepoint and the fact that the radius of injection depends on where we are in $X$.  However, in the compact setting, we have a uniform injectivity radius, so we may may avoid this step altogether.
Morally, uniformity in the basepoint is due to the fact that in the compact setting, there are no proper invariant subspaces near which an orbit can become trapped for long periods of time.
\end{rem}


\section{Equidistribution for Arithmetic Sequences Along Abelian Horospherical Flows}\label{sect:EquiArith}

Let $G = \slnr$, $\Gamma = \slnz$, and $X = \Gamma\backslash G$. Let $U$ be an upper triangular unipotent subgroup of the form
\begin{align}
U=
\left\{\left(\begin{matrix}
I_{m} & \vline & \bigast \\
\hline
0 & \vline & I_{n-m}
\end{matrix} 
\right)\right\}\label{eq:U}
\end{align}
for $m<n$.  Note that $U\cong \R^d$ as groups for $d=m(n-m)$ under any identification $u(\bold{t})$ of ${\bf t}\in\R^d$ with the upper-right block.  Recall that the Haar measure on $U$ is the Lebesgue measure on $\R^d$ under this identification, which we normalize so that $u([0,1]^d)$ has unit measure.  Observe that $U$ is horospherical with respect to the element 
\begin{align*}
a_t = {\rm  diag}(\underbrace{e^{t(n-m)/n}, \cdots, e^{t(n-m)/n}}_{m}, \underbrace{e^{-tm/n}, \cdots, e^{-tm/n}}_{n-m})
\end{align*}
for any $t>0$ and that conjugation by $a_t$ scales all entries in the upper-right block of $U$ by $e^{t(n-m)/n}e^{tm/n}=e^t$.  Hence, for this choice of $a_t$, we have $B_T = a_{\log T} u([0,1]^d) a_{-\log T} = u([0,T]^d)$.  For this reason we will conflate the notation and write $B_T$ for both $[0,T]^d\subseteq\R^d$ and $u([0,T]^d)\subseteq U$.

Let $\psi$ be an additive character of $U$ (so $\psi(\bold{t}) = e^{i{\bf a}\cdot \bf{t}}$ for some ${\bf a}\in\R^d$).  Define measure $\nu_T$ and (complex) measure $\mu_{T,\psi}$ on $X$ via duality: for $f\in C_c^\infty(X)$ let
\begin{align*}
\int_X f d\nu_T = \nu_T(f) := \frac{1}{|B_T|}\int_{B_T} f(x_0 u({\bf t})) d{\bf t}
\end{align*}
and
\begin{align*}
\int_X f d\mu_{T,\psi}=\mu_{T,\psi}(f) := \frac{1}{|B_T|}\int_{B_T} \psi({\bf t})\left(f(x_0 u({\bf t}))-\int_X f dm_X\right) d{\bf t}.
\end{align*}

Our main goal in this section is to obtain an effective rate of equidistribution along (multivariate) arithmetic sequences of inputs for the right action of $U$ on $X$.  To do this, we first present the following lemma, the proof of which closely follows the proof of Lemma 3.1 in \cite{VenkSparse} for the case of $G=\slr$ and $\Gamma$ cocompact.

\begin{lem}\label{lem:VenkLem}
Let $x_0 =\Gamma g_0\in X $ satisfy (\ref{eq:dio1}) for $T>R>C$.  Then there exists $b>0$
 such that for any $f\in C^\infty_c(X)$ and additive character $\psi$, 
\begin{align*}
\left|\mu_{T,\psi}(f)\right|\ll R^{- b} \sob{\infty}{\ell}(f)
\end{align*}
where $\ell$ is as in Theorem \ref{thm:equidist}.
\end{lem}

\begin{rmk}
As noted in \cite{VenkSparse}, the significance of this lemma is that the implicit constant is independent of choice of $\psi$.  This can be shown for highly oscillatory $\psi$ using integration by parts and for almost constant $\psi$ using equidistribution of the horospherical flow directly, thus this lemma is most significant for $\psi$ of moderate oscillation.
The proof will use our effective equidistribution result as well as a variety of technical integral manipulations that nonetheless do not require any heavy machinery.
\end{rmk}

\begin{proof}

Let $1\leq H\leq T$ and define a complex measure $\sigma_H$ on $U$ by
\begin{align*}
\int_U g d\sigma_H = \sigma_H(g) := \frac{1}{|B_H|}\int_{B_H} \overline\psi({\bf t})g(u({\bf t}))d{\bf t}
\end{align*}
for $g \in C_c^\infty(U)$.

Let $f\ast\sigma_H$ be the right convolution of $f$ by $\sigma_H$, i.e., for $x\in X$
\begin{align*}
f\ast \sigma_H (x) &= \int f(x u({\bf t})^{-1})d\sigma_H({\bf t})\\
&=\frac{1}{|B_H|} \int_{B_H} \overline\psi ({\bf t})f(x u({\bf t})^{-1}) d{\bf t}.
\end{align*}
Notice that by switching the order of integration (one may verify that the conditions of Fubini's theorem are satisfied) and using invariance of the Haar measure, we have
\begin{align*}
\int_X f\ast\sigma_H dm_X &= \int_X \frac{1}{|B_H|} \int_{B_H} \overline\psi ({\bf t})f(x u({\bf t})^{-1}) d{\bf t} \hspace{.03in}dm_X(x)\\
&= \frac{1}{|B_H|} \int_{B_H} \overline\psi ({\bf t})\left(\int_Xf(x u({\bf t})^{-1}) dm_X(x)\right)d{\bf t}\\
&= \frac{1}{|B_H|} \int_{B_H} \overline\psi ({\bf t})\left(\int_X f dm_X \right)d{\bf t}.
\end{align*}
Hence,
\begin{align*}
\mu_{T,\psi}(f\ast \sigma_H) &= \frac{1}{|B_T|}\int_{B_T}\psi({\bf t})\left( f\ast\sigma_H(x_0u({\bf t})) - \int_X f\ast\sigma_H dm_X \right)d{\bf t}\\
&= \frac{1}{|B_T|}\int_{B_T}\psi({\bf t})\frac{1}{|B_H|}\int_{B_H} \overline\psi({\bf s})\left(f(x_0 u({\bf t})u({\bf s})^{-1})-\int_X f dm_X\right)d{\bf s}d{\bf t}\\
&= \frac{1}{|B_T||B_H|}\int_{B_T}\int_{B_H}\psi({\bf t-s})\left(f(x_0 u({\bf t-s}))-\int_X f dm_X\right)d{\bf s}d{\bf t}
\end{align*}
since $\overline\psi ({\bf s}) =\psi (-{\bf s})$ and $U\cong R^d$.  Now by switching the order of integration and applying a change of variables, we get
\begin{align*}
\mu_{T,\psi}(f\ast \sigma_H)  &= \frac{1}{|B_T||B_H|}\int_{B_H}\int_{B_T- {\bf s}}\psi({\bf t})\left(f(x_0 u({\bf t}))-\int_X f dm_X\right)d{\bf t}d{\bf s}.
\end{align*}
But we may also write
\begin{align*}
\mu_{T,\psi}(f) &= \frac{1}{|B_T|}\int_{B_T} \psi({\bf t})\left(f(x_0 u({\bf t}))-\int_X f dm_X\right) d{\bf t}\\
&=\frac{1}{|B_T||B_H|}\int_{B_H}\int_{B_T} \psi({\bf t})\left(f(x_0 u({\bf t}))-\int_X f dm_X\right) d{\bf t}d{\bf s}.
\end{align*}
Thus
\begin{align*}
&|\mu_{T,\psi}(f)-\mu_{T,\psi}(f\ast \sigma_H)|\\
&\leq \frac{1}{|B_T||B_H|}\int_{B_H}\int_{B_T\triangle (B_T-{\bf s})} \left| f(x_0 u({\bf t}))-\int_X f dm_X\right| d{\bf t}d{\bf s}\\
&\ll \frac{1}{|B_T||B_H|}\int_{B_H} |B_T\triangle (B_T-{\bf s})| \sob{\infty}{0}(f) d{\bf s}.
\end{align*}
But notice that $B_T\triangle (B_T-{\bf s})$ is simply the symmetric difference of two shifted cubes, the measure of which will be maximized when ${\bf s} = (H, \cdots, H)$ (see Figure \ref{fig:symmdiff}). Hence, 
\begin{align*}
|B_T\triangle (B_T-{\bf s})| &\leq 2 (T^d - (T-H)^d)\\
&= 2( d T^{d-1}H - \cdots \pm d T H^{d-1} \mp H^d)\\
&\ll T^{d-1}H.
\end{align*}
since $H\leq T$ implies that the leading term dominates.

\begin{figure}
\begin{center}
\begin{tikzpicture}[scale=.8]
\draw [-] (0,0) -- (5,0);
\draw [-] (0,0) -- (0,5);
\draw [-] (0,5) -- (5,5);
\draw [-] (5,0) -- (5,5);
\draw [-] (1,1) -- (6,1);
\draw [-] (1,1) -- (1,6);
\draw [-] (1,6) -- (6,6);
\draw [-] (6,1) -- (6,6);
\draw[fill=gray!20]  (0,0) -- (0,5) -- (1,5) -- (1,1) -- (5,1) -- (5,0) -- cycle;
\draw[fill=gray!20]  (6,6) -- (6,1) -- (5,1) -- (5,5) -- (1,5) -- (1,6) -- cycle;
\draw [red,decorate,decoration={brace,amplitude=10pt,mirror}]
(0,0) -- (5,0) node [red,midway,below,yshift=-10pt] 
{$T$};
\draw [red,decorate,decoration={brace,amplitude=10pt}]
(0,0) -- (0,5) node [red,midway,left,xshift=-10pt] 
{$T$};
\draw [red,decorate,decoration={brace,amplitude=5pt,mirror}]
(5,1) -- (6,1) node [red,midway,below,yshift=-5pt] 
{$H$};
\draw [red,decorate,decoration={brace,amplitude=5pt}]
(1,5) -- (1,6) node [red,midway,left,xshift=-5pt] 
{$H$};
\draw [red,decorate,decoration={brace,amplitude=7pt,mirror}]
(1,1) -- (1,5) node [red,midway,right,xshift=7pt] 
{$T-H$};
\draw [red,decorate,decoration={brace,amplitude=7pt}]
(1,1) -- (5,1) node [red,midway,above,yshift=5pt] 
{$T-H$};
\end{tikzpicture}
\caption{The symmetric difference between $B_T$ and $B_T - (H,\cdots,H)$.}
\label{fig:symmdiff}
\end{center}
\end{figure}
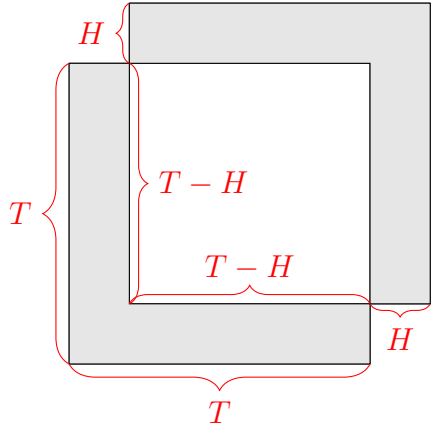

It follows that $\int_{B_H} |B_T\triangle (B_T-{\bf s})| d{\bf s} \ll T^{d-1} H^{d+1}$.  Thus, 
\begin{align}
|\mu_{T,\psi}(f)-\mu_{T,\psi}(f\ast \sigma_H)| \ll \frac{T^{d-1}H^{d+1}}{|B_T||B_H|} \sob{\infty}{0}(f) = \frac{H}{T}\sob{\infty}{0}(f).\label{eq:H/T}
\end{align}

Now consider 
\begin{align*}
|\mu_{T,\psi}(f\ast \sigma_H)|^2 &= \left|\frac{1}{|B_T|}\int_{B_T} \psi({\bf t})\left(f\ast\sigma_H(x_0u({\bf t}))-\int_X f\ast \sigma_H dm_X\right)d{\bf t}\right|^2\\
&\leq \frac{1}{|B_T|^2}\left(\int_{B_T} \left|f\ast\sigma_H(x_0u({\bf t}))-\int_X f\ast \sigma_H dm_X\right|d{\bf t}\right)^2\\
&= \frac{1}{|B_T|^2}\left<1, \left|f\ast\sigma_H(x_0u(\cdot))-\int_X f\ast \sigma_H dm_X\right| \right>_{L^2(B_T)}^2.
\end{align*}
By Cauchy-Schwarz, we know that
\begin{align*}
|\mu_{T,\psi}(f\ast \sigma_H)|^2 &\leq \frac{1}{|B_T|^2} \nn{1}_{L^2(B_T)}^2 \nn{f\ast\sigma_H(x_0u(\cdot))-\int_X f\ast \sigma_H dm_X}_{L^2(B_T)}^2.
\end{align*}
Now, $\nn{1}_{L^2(B_T)}^2 = \int_{B_T} 1^2 d{\bf t} = |B_T|$, and 
\begin{align*}
\nn{f\ast\sigma_H(x_0u(\cdot))-\int_X f\ast \sigma_H dm_X}_{L^2(B_T)}^2 &= \int_{B_T} \left|f\ast\sigma_H(x_0u({\bf t}))-\int_X f\ast \sigma_H dm_X\right|^2 d{\bf t}\\
&=|B_T|\nu_T\left(\left|f\ast\sigma_H-\int_X f\ast \sigma_H dm_X\right|^2\right)
\end{align*}
which shows that
\begin{align}
|\mu_{T,\psi}(f\ast \sigma_H)|^2 \leq \nu_T\left(\left|f\ast\sigma_H-\int_X f\ast \sigma_H dm_X\right|^2\right).\label{eq:munu}
\end{align}
Hence, by (\ref{eq:H/T}) and (\ref{eq:munu}), we have
\begin{align}
|\mu_{T,\psi}(f)| &\leq |\mu_{T,\psi}(f)-\mu_{T,\psi}(f\ast \sigma_H)| + |\mu_{T,\psi}(f\ast \sigma_H)|\nonumber\\
&\ll \frac{H}{T}\sob{\infty}{0}(f) + \nu_T\left(\left|f\ast\sigma_H-\int_X f\ast \sigma_H dm_X\right|^2\right)^{1/2}. \label{eq:muTpsi}
\end{align}

To estimate $\nu_T\left(\left|f\ast\sigma_H-\int_X f\ast \sigma_H dm_X\right|^2\right)$, observe that
\begin{align*}
&\left|f\ast\sigma_H(x)-\int_X f\ast \sigma_H dm_X\right|^2\\
&= \left|\frac{1}{|B_H|}\int_{B_H} \overline\psi({\bf s})\left( [u({\bf s})f](x)-\int_X f dm_X \right) d{\bf s}\right|^2\\
&=\frac{1}{|B_H|^2}\left(\int_{B_H} \overline\psi({\bf s_1})\left([u({\bf s_1})f](x)-\int_X f dm_X \right)d{\bf s_1}\right)\\
&\hspace{.55in}\cdot\left(\int_{B_H} \psi({\bf s_2})\left([u({\bf s_2})\overline f](x)-\int_X \overline f dm_X\right)  d{\bf s_2}\right)\\
&= \frac{1}{|B_H|^2}\int_{B_H}\int_{B_H}\psi({\bf s_2-s_1}) \left[\left([u({\bf s_1})f](x)-\int_X f dm_X \right)\hspace{-4pt}\left([u({\bf s_2})\overline f](x)-\int_X \overline f dm_X\right)\right]d{\bf s_1}d{\bf s_2}.
\end{align*}
When we apply $\nu_T$ to this, we can change the order of integration so that the innermost integral is over $B_T$, with the character $\psi({\bf s_2-s_1})$ outside this integral.  We may then integrate separately over the four terms we get by expanding the bracketed product above.  That is, 
\begin{align}
&\nu_T\left(\left|f\ast\sigma_H-\int_X f\ast \sigma_H dm_X\right|^2\right)\nonumber\\
&=\frac{1}{|B_H|^2}\int_{B_H}\int_{B_H}\psi({\bf s_2-s_1}) \nu_T\left(\left(u({\bf s_1})f-\int_X f dm_X \right)\left(u({\bf s_2})\overline f-\int_X \overline f dm_X\right)\right)d{\bf s_1}d{\bf s_2}\nonumber\\
&=\frac{1}{|B_H|^2}\int_{B_H}\int_{B_H}\psi({\bf s_2-s_1}) F({\bf s_1, s_2}) d{\bf s_1}d{\bf s_2}\label{eq:introF}
\end{align}
where
\begin{align}
F({\bf s_1, s_2}) &= \nu_T(u({\bf s_1})f\cdot u({\bf s_2})\overline f)\nonumber\\
&-  \nu_T(u({\bf s_1})f) \int_X \overline f dm_X \nonumber\\
&-  \nu_T(u({\bf s_2})\overline f) \int_X  f dm_X\label{eq:F}\\
&+ \left|\int_X f dm_X\right|^2.\nonumber
\end{align}

Now from Theorem \ref{thm:equidist} we know that for arbitrary $\tilde f\in C_c^\infty(X)$ and $x_0$ satisfying the Diophantine basepoint condition (\ref{eq:dio1}) with $T>R>C$, we have
\begin{align*}
\left|\nu_T(\tilde f)-\int_X \tilde f dm_X\right| = \left|\frac{1}{|B_T|}\int_{B_T} \tilde f(x_0 u({\bf t})) d{\bf t}-\int_X \tilde f dm_X\right| \ll  R^{-\gamma}\sob{\infty}{\ell}(\tilde f),
\end{align*}
that is,
\begin{align}
\nu_T(\tilde f) = \int_X \tilde f dm_X +\mathcal{O}(R^{-\gamma}\sob{\infty}{\ell}(\tilde f)).\label{eq:equidist}
\end{align}
Applying this to the function $\tilde f = u({\bf s_1})f$, we find that
\begin{align*}
\nu_T(u({\bf s_1})f) = \int_X u({\bf s_1})f dm_X + \mathcal{O}(R^{-\gamma}\sob{\infty}{\ell}(u({\bf s_1})f)).
\end{align*}
But since $m_X$ is the Haar measure, 
\begin{align*}
\int_X u({\bf s_1})f dm_X = \int_X f(xu({\bf s_1})^{-1}) dm_X(x) = \int_X f dm_X.
\end{align*}
Thus
\begin{align*}
\nu_T(u({\bf s_1})f) \int_X \overline f dm_X = \left|\int_X f dm_X\right|^2 + \mathcal{O}\left(R^{-\gamma}\sob{\infty}{\ell}(u({\bf s_1})f)\left|\int_X \overline f dm_X\right|\right).
\end{align*}
Furthermore, from Sobolev norm property (\ref{Sob3}), we know that for $f\in C^\infty_c(X)$ and $h\in G$, we have $\sob{\infty}{\ell}(hf) \ll_\ell \nn{h}^\ell \sob{\infty}{\ell}(f)$, where $\nn{h}$ is the operator norm of ${\rm  Ad}_{h^{-1}}$.  Since the entries of $u({\bf s})^{-1}$ are bounded by $\max(1,|{\bf s}|)$, we have $\nn{u({\bf s})} \ll \max(1,|{\bf s}|)^2$.  Thus for ${\bf s_1} \in [0,H]$ with $H\geq 1$, $\sob{\infty}{\ell}(u({\bf s_1})f) \ll H^{2\ell}\sob{\infty}{\ell}(f)$.  Combining this with the bound $\left|\int \overline f dm_X\right| \ll \sob{\infty}{0}(f) \ll \sob{\infty}{\ell}(f)$, we find that
\begin{align*}
\nu_T(u({\bf s_1})f) \int_X \overline f dm_X = \left|\int_X f dm_X\right|^2 + \mathcal{O}(R^{-\gamma}H^{2\ell}\sob{\infty}{\ell}(f)^2).
\end{align*}
Likewise,
\begin{align*}
\nu_T(u({\bf s_2})\overline f) \int_X f dm_X = \left|\int_X f dm_X\right|^2 + \mathcal{O}(R^{-\gamma}H^{2\ell}\sob{\infty}{\ell}(f)^2).
\end{align*}
Therefore, (\ref{eq:F}) becomes simply
\begin{align*}
F({\bf s_1, s_2}) =\nu_T(u({\bf s_1})f\cdot u({\bf s_2})\overline f)-\left|\int_X f dm_X\right|^2 + \mathcal{O}(T^{-\alpha\gamma}H^{2\ell}\sob{\infty}{\ell}(f)^2).
\end{align*}
Substituting this back into (\ref{eq:introF}), we conclude that
\begin{align}
&\nu_T\left(\left|f\ast\sigma_H(x)-\int_X f\ast \sigma_H dm_X\right|^2\right)\nonumber\\
&\ll \frac{1}{|B_H|^2}\int_{B_H}\int_{B_H}\left|\nu_T(u({\bf s_1})f\cdot u({\bf s_2})\overline f)-\left|\int_X f dm_X\right|^2 \right|d{\bf s_1}d{\bf s_2}+R^{-\gamma}H^{2\ell}\sob{\infty}{\ell}(f)^2.\label{eq:backtointroF}
\end{align}
But now notice that
\begin{align*}
\int_X u({\bf s_1})f\cdot u({\bf s_2})\overline f dm_X = \left<u({\bf s_1})f,u({\bf s_1})f\right>_{L^2(X)} = \left<u({\bf s_1-s_2})f,f\right>_{L^2(X)}
\end{align*}
so by the triangle inequality, we can estimate
\begin{align}
\left|\nu_T(u({\bf s_1})f\cdot u({\bf s_2})\overline f)-\left|\int_X f dm_X\right|^2 \right| \leq& \left|\nu_T(u({\bf s_1})f\cdot u({\bf s_2})\overline f)-\int_X u({\bf s_1})f\cdot u({\bf s_2})\overline f dm_X  \right|\nonumber\\
&+ \left|\left<u({\bf s_1-s_2})f,f\right>_{L^2(X)}-\left|\int_X f dm_X\right|^2\right|.\label{eq:longtriangleineq}
\end{align}
Again, by our equidistribution result in (\ref{eq:equidist}), we know that 
\begin{align}
\left|\nu_T(u({\bf s_1})f\cdot u({\bf s_2})\overline f)-\int_X u({\bf s_1})f\cdot u({\bf s_2})\overline f dm_X  \right|\ll R^{-\gamma}\sob{\infty}{\ell}(u({\bf s_1})f\cdot u({\bf s_2})\overline f)\label{eq:oddequidist}
\end{align}
and by properties (\ref{Sob2}) and (\ref{Sob3}) of Sobolev norms, we have
\begin{align}
\sob{\infty}{\ell}(u({\bf s_1})f\cdot u({\bf s_2})\overline f) \ll \sob{\infty}{\ell}(u({\bf s_1})f) \sob{\infty}{\ell}(u({\bf s_2})f) \ll H^{4\ell} \sob{\infty}{\ell}(f)^2\label{eq:Sobbd}
\end{align}
 for ${\bf s_1, s_2}\in [0,H]$.  Thus, from (\ref{eq:Sobbd}), (\ref{eq:oddequidist}), and (\ref{eq:longtriangleineq}), equation (\ref{eq:backtointroF}) becomes
\begin{align}
&\nu_T\left(\left|f\ast\sigma_H(x)-\int_X f\ast \sigma_H dm_X\right|^2\right) \nonumber\\ 
&\ll \frac{1}{|B_H|^2}\int_{B_H}\int_{B_H}\left|\left<u({\bf s_1-s_2})f,f\right>_{L^2(X)}-\left|\int_X f dm_X\right|^2\right|d{\bf s_1}d{\bf s_2}+ R^{-\gamma}H^{4\ell}\sob{\infty}{\ell}(f)^2.\label{eq:lastest}
\end{align}
Now from Corollary \ref{cor:expmixing} (\ref{umtxcoeff}), we know there exists $\beta>0$ such that for any ${\bf s}\in \R^d$,
\begin{align}
\left|\left<u({\bf s})f,f\right>_{L^2(X)}-\left|\int_X f dm_X\right|^2\right| \ll \max(1,|{\bf s}|)^{-\beta}\sob{\infty}{\ell}(f)^2.\label{eq:spectralGap}
\end{align}
Then for ${\bf s} = {\bf s_1}-{\bf s_2}$, we have the following problem: We want to bound the integral in (\ref{eq:lastest}) by a power of $H$, but for $({\bf s_1},{\bf s_2})$ close to the diagonal in $B_H\times B_H$ we cannot do better than a constant times the Sobolev norm of $f$ in (\ref{eq:spectralGap}).  We will address this by integrating separately over a neighborhood of the diagonal that has small measure (depending on $H$) and away from the diagonal where $\max(1,|{\bf s_1} - {\bf s_2}|)$ is dominated by $H$.

To make this precise, let $D:=\{  ({\bf s_1,s_2})\in B_H \times B_H \hspace{2pt}\vert\hspace{3pt} {\bf s_1 = s_2} \}$ be the diagonal of $ B_H \times B_H $ and define $D_\epsilon:=\{ ({\bf s_1,s_2})\in B_H\times B_H \hspace{2pt}\vert\hspace{3pt} |{\bf s_1-s_2| }<\epsilon \}$.
Notice that $D$ is a $d$-dimensional subset of $\R^{2d}$ with diameter $\sqrt{2d}H$.  Furthermore, any point satisfying $|{\bf s_1 - s_2}| = \epsilon$ is distance $\epsilon/\sqrt{2}$ from the diagonal, so $D_\epsilon$ is an ($\epsilon/\sqrt{2}$)-neighborhood of $D$ sitting inside $[0,H]^{2d}$.  Thus $D_\epsilon$ is contained within a box in $\R^{2d}$ with $d$ side-lengths of $\sqrt{2d}H$ and $d$ side-lengths of $2\epsilon/\sqrt{2}$, so
\begin{align*}
|D_\epsilon| \ll H^d\epsilon^d
\end{align*}
(see Figure \ref{fig:nbhdofdiag}). In particular, if $\epsilon = H^\zeta$ (for $0<\zeta<1$ to be determined), then 
\begin{align*}
\left|\{ ({\bf s_1}, {\bf s_2})\in B_H\times B_H \hspace{2pt}\vert\hspace{3pt} |{\bf s_1-s_2| }<H^\zeta \}\right| \ll H^{d(1+\zeta)}.
\end{align*}
In this region, the integrand is dominated by 1, so when we integrate over this region and divide by $|B_H|^2 = H^{2d}$ (as we are doing in  (\ref{eq:lastest})), we get a term of order $H^d({\zeta-1)}\sob{\infty}{\ell}(f)^2$.  

\begin{figure}
\begin{center}
\begin{tikzpicture}
\draw[fill=gray!20]  (0,0) -- (1,0) -- (5,4) -- (5,5) -- (4,5) -- (0,1) -- cycle;
\draw [-] (0,0) -- (5,0);
\draw [-] (0,0) -- (0,5);
\draw [-] (0,5) -- (5,5);
\draw [-] (5,0) -- (5,5);
\draw [-,red] (0,0) -- (5,5);
\draw [-] (-.5,.5) -- (.5,-.5);
\draw [-] (4.5,5.5) -- (5.5,4.5);
\draw [-] (-.5,.5) -- (4.5,5.5);
\draw [-] (.5,-.5) -- (5.5,4.5);
\draw [-, red] (2.5,2.5) -- (3,2);
\draw [-, dashed,gray] (2,2) -- (2,0);
\draw [-, dashed,gray] (3,2) -- (3,0);
\draw [-, dashed,gray] (3,2) -- (0,2);
\draw [red,decorate,decoration={brace,amplitude=10pt}]
(-.5,.5,0) -- (4.5,5.5) node [red,midway,xshift=-17pt,yshift=15pt] 
{$\sqrt{2d}H$};
\draw [decorate,decoration={brace,amplitude=5pt,mirror}]
(2,0) -- (3,0);
\node[below] at (4.9,0) {$H$};
\node[left] at (0,4.9) {$H$};
\node[right,red] at (2.5,2.55) {\small $\epsilon/\sqrt{2}$};
\node[below] at (2,0) {\small ${\bf s_2}$};
\node[left] at (0,2) {\small ${\bf s_2}$};
\node[below] at (3.1,0) {\small ${\bf s_1}$};
\node[below] at (2.5,-.1) {\small $\epsilon$};
\node[gray] at (4,4.5) {$D_\epsilon$};
\node[red] at (4,3.6) {$D$};
\end{tikzpicture}
\caption{The measure of the set where $|{\bf s_1-s_2| }<\epsilon$ has measure bounded by $H^d\epsilon^d$ in $B_H\times B_H$ (shown here for one dimensional $U$).}
\label{fig:nbhdofdiag}
\end{center}
\end{figure}
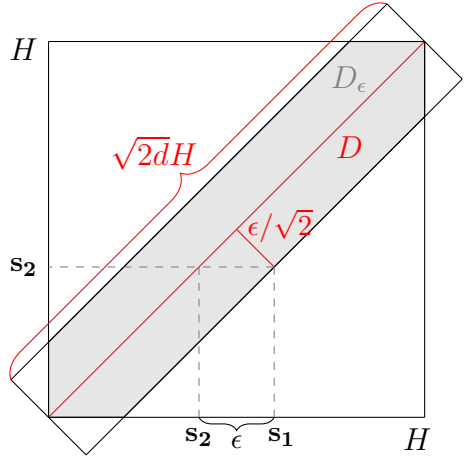

On the other hand, for $|{\bf s_1- s_2}|\geq H^\zeta$, 
we can say that
\begin{align}
\left|\left<u({\bf s_1-s_2})f,f\right>_{L^2(X)}-\left|\int_X f dm_X\right|^2\right| &\ll \max(1,|{\bf s_1- s_2}|)^{-\beta}\sob{\infty}{\ell}(f)^2\nonumber\\
&\leq H^{-\zeta\beta}\sob{\infty}{\ell}(f)^2.\nonumber
\end{align}
Hence,
\begin{align}
 \frac{1}{|B_H|^2}\int_{B_H}\int_{B_H}\left|\left<u({\bf s_1-s_2})f,f\right>_{L^2(X)}-\left|\int_X f dm_X\right|^2\right| d{\bf s_1}d{\bf s_2} &\ll (H^{-\zeta\beta} + H^{d(\zeta-1)})\sob{\infty}{\ell}(f)^2\nonumber\\
&= H^{-d\beta/(d+\beta)}\sob{\infty}{\ell}(f)^2\label{eq:betaBd}
\end{align}
where we have chosen $\zeta=d/(d+\beta)$ to optimize the error.

Together, the bounds in (\ref{eq:lastest}) and (\ref{eq:betaBd}) imply that
\begin{align}
\nu_T\left(\left|f\ast\sigma_H(x)-\int_X f\ast \sigma_H dm_X\right|^2\right) \ll& \left(R^{-\gamma}H^{4\ell} + H^{-d\beta/(d+\beta)}\right)\sob{\infty}{\ell}(f)^2.\label{eq:nuTbd}
\end{align}
Finally, from (\ref{eq:muTpsi}) and (\ref{eq:nuTbd}), we have
\begin{align*}
|\mu_{T,\psi}(f)| \ll \left( T^{-1}H + R^{-\gamma/2}H^{2\ell} + H^{-d\beta/(2d+2\beta)}\right)\sob{\infty}{\ell}(f).
\end{align*}
Since $\gamma<1$ and $R<T$, 
the first term decays more quickly that the second, and can be ignored.  Thus the decay is optimized when
\begin{align*}
H^{-d\beta/(2d+2\beta)} &= R^{-\gamma/2} H^{2\ell}\\
H &= R^{\gamma(d+\beta)/(4\ell d + 4\ell\beta + d\beta)}.
\end{align*}
This demonstrates the claim that
\begin{align*}
|\mu_{T,\psi}(f)| &\ll R^{-b} \sob{\infty}{\ell}(f)
\end{align*}
where $b = d\beta\gamma/(8d\ell +8\ell\beta+2d\beta)$.
\end{proof}


We will now use this lemma to establish an effective equidistribution bound along multivariate arithmetic sequences.

Let $K_1, \dots, K_d\geq1$ and define $K$ to be the diagonal matrix
\begin{align*}
K := {\rm  diag}(K_1,\cdots, K_d) = \begin{pmatrix}K_1&&\\&\ddots&\\&&K_d\end{pmatrix}
\end{align*}
and $|K| = \det(K) = K_1K_2\cdots K_d$.

We want to understand the behavior of
\begin{align}
S:=\sum_{\substack{{\bf k}\in\Z^d\\ K{\bf k}\in B_T}} f(x_0 u(K{\bf k})).\label{eq:Sdef}
\end{align}
For equidistribution, we want this to be close to $\#\{{\bf k}\in\Z^d | K{\bf k}\in B_T\} \int_X f dm_X \approx \frac{T^d}{|K|}\int_X f dm_X$.  For $x_0$ satisfying a basepoint property, we have the following result.

\begin{thm}
Let $K = {\rm  diag}(K_1,\cdots, K_d)$ with $T\geq K_1, \dots, K_d\geq1$ and determinant $|K|$.  Then for all $x_0\in X$ satisfying (\ref{eq:dio1}) with $T>R>C_K$, we have
\begin{align*}
\left|\sum_{\substack{{\bf k}\in\Z^d\\ K{\bf k}\in B_T}} f(x_0 u(K{\bf k}))-\frac{T^d}{|K|}\int_X f dm_X\right| \ll \left(T^d R^{-b/(d+1)} |K|^{-d/(d+1)} +\frac{T^{d-1}\max_i K_i}{|K|}\right)\sob{\infty}{\ell}(f)
\end{align*}
where $C_K:=\max(C, (2/\min K_i)^{(d+1)/b}|K|^{1/b})$ with $C$ and $\ell$ as in Theorem \ref{thm:equidist}.\label{thm:arithequidist}
\end{thm}

\begin{proof}

Let $\delta>0$ be small (to be determined) and define the single-variable hat function
\begin{align*}
g_\delta(t) := \max(\delta^{-2}(\delta-|t|), 0)
\end{align*}
for $t\in\R$ and (through slight abuse of notation) the multivariable function
\begin{align*}
g_\delta({\bf t}) := g_\delta(t_1)\cdots g_\delta(t_d)
\end{align*}
for ${\bf t} = (t_1, \cdots, t_d) \in \R^d$.  Notice that $\int_{\R^d} g_\delta ({\bf t}) d{\bf t} = 1$ and $\text{supp } (g_\delta) \subseteq [-\delta, \delta]^d$.

Define an approximation to the sum $S$ by
\begin{align}
S_{{\rm  approx}} := \int_{B_T} \left(\sum_{{\bf k}\in\Z^d} g_\delta({\bf t}-K{\bf k})\right) f(x_0 u({\bf t}))d{\bf t}.\label{eq:Sapproxdef}
\end{align}
That is, instead of averaging $f$ over the lattice points of $K\Z^d$, we average over small neighborhoods around the lattice points using the bump function $g_\delta$, since $\sum_{{\bf k}} g_\delta({\bf t}-K{\bf k})$ is supported on a disjoint union of $\delta$-cubes centered around the points of $K\Z^d$ (that is, so long as $\delta < \min_i K_i/2$).

We want to show that $S_{\rm  approx}$ can be written
\begin{align}
 S_{\rm  approx} = \left( \sum_{\substack{{\bf k}\in\Z^d\\ K{\bf k}\in B_T}} \int_{[-\delta, \delta]^d + K{\bf k}} g_\delta({\bf t}-K{\bf k}) f(x_0u({\bf t}))d{\bf t}\right) +  r(T,K,f,d)\label{eq:Sapproxapprox}
\end{align}
where $r(T,K,f,d)$ is an error term depending on $T$, $K_1, \cdots, K_d$, $f$, and dimension $d$.  To see this, observe that in both (\ref{eq:Sapproxdef}) and (\ref{eq:Sapproxapprox}) we are integrating $f$ against a sum of bump functions supported on a disjoint union of $\delta$-cubes centered at the lattice points of $K\Z^d$.  However, in (\ref{eq:Sapproxdef}) we are integrating over the region shaded in red Figure \ref{fig:suminterror}, whereas in (\ref{eq:Sapproxapprox}) we are integrating over the region shaded in gray (that is, we are only integrating against the bump functions whose centers intersect $B_T$).

\begin{figure}
\begin{center}
\begin{tikzpicture}
\foreach \Point in {(0,0),(2.9,0),(5.8,0),
(0,2.5),(2.9,2.5), (5.8,2.5),
(0,5), (2.9,5),(5.8,5),,
(0,7.5),(2.9,7.5),(5.8,7.5)}
{
	\begin{scope}
  		\clip (0,0) rectangle (8,8);
  		\draw[pattern=north east lines, pattern color=red!30] \Point +(-1,-1) rectangle +(1,1);
	\end{scope}
	\draw[pattern=north west lines, pattern color=gray!60] \Point +(-1,-1) rectangle +(1,1);
	\draw \Point +(-1,-1) rectangle +(1,1) ;
	\fill \Point circle[radius=2pt];
}
\foreach \Point in {(8.7,0),(8.7,2.5),(8.7,5),(8.7,7.5)}
{
	\begin{scope}
  		\clip (0,0) rectangle (8,8);
  		\draw[pattern=north east lines, pattern color=red!30] \Point +(-1,-1) rectangle +(1,1);
	\end{scope}
	\fill \Point circle[radius=2pt];
	\draw \Point +(-1,-1) rectangle +(1,1) ;
}
\draw [-] (0,0) -- (8,0);
\draw [-] (0,0) -- (0,8);
\draw [-] (0,8) -- (8,8);
\draw [-] (8,0) -- (8,8);
\draw [red,decorate,decoration={brace,amplitude=10pt,mirror}]
(0,0) -- (2.9,0) node [midway,yshift=-17pt] 
{$K_1$};
\draw [red,decorate,decoration={brace,amplitude=10pt}]
(0,0) -- (0,2.5) node [midway,xshift=-18pt] 
{$K_2$};
\draw [red,decorate,decoration={brace,amplitude=7pt,mirror}]
(2.9,0) -- (3.9,0) node [midway,yshift=-13pt] 
{$\delta$};
\draw [red,decorate,decoration={brace,amplitude=7pt}]
(0,2.5) -- (0,3.5) node [midway,xshift=-13pt] 
{$\delta$};
\draw [red,decorate,decoration={brace,amplitude=10pt,aspect=.55}]
(0,8) -- (8,8) node [midway,yshift=16pt,xshift=12pt] 
{$T$};
\draw [red,decorate,decoration={brace,mirror,amplitude=10pt,aspect=.47}]
(8,0) -- (8,8) node [midway,yshift=-7pt,xshift=17pt] 
{$T$};
\end{tikzpicture}
\caption{The area shaded in red indicates the region over which we are integrating in the definition of $S_{\rm  approx}$, whereas the area shaded in gray represents the region over which we are integrating in our estimate of $S_{\rm  approx}$ given in (\ref{eq:Sapproxapprox}).  The difference between the two integrals can be bounded by the number of $\delta$-cubes intersecting the boundary of $B_T$ multiplied by the supremum of $f$.}
\label{fig:suminterror}
\end{center}
\end{figure}

Thus all of the possible error comes from integrating over those $\delta$-cubes that intersect the boundary of $B_T$.  Consider a face of $B_T$ that is orthogonal to the $i^{\rm  th}$ standard basis vector. It will intersect at most $T/K_j + \mathcal{O}(1)$ of these cubes along an edge in the $j^{\rm  th}$ direction for $j\neq i$.  Hence, the total number of cubes that face intersects can be bounded by 
\begin{align*}
\frac{T}{K_1}\cdots\frac{T}{K_{i-1}}\cdot\frac{T}{K_{i+1}}\cdots\frac{T}{K_d} = T^{d-1}\frac{K_i}{|K|}.
\end{align*}
Since $g_\delta$ integrates to one, the error that results from integrating over one of these $\delta$-cubes is bounded by $\sob{\infty}{0}(f)$.  Then considering all the faces of $B_T$, we see that the error satisfies
\begin{align*}
|r(T,K,f,d)| \ll \sob{\infty}{0}(f)\sum_{i=1}^d T^{d-1}\frac{K_i}{|K|} \ll T^{d-1}\frac{\max_i K_i}{|K|}\sob{\infty}{0}(f).
\end{align*}
Then by a change of variables in (\ref{eq:Sapproxapprox}), we have
\begin{align}
S_{{\rm  approx}} &= \left( \sum_{\substack{{\bf k}\in\Z^d\\ K{\bf k}\in B_T}} \int_{[-\delta, \delta]^d} g_\delta({\bf s}) f(x_0u(K{\bf k}+{\bf s}))d{\bf s}\right) +  r(T,K,f,d). \label{eq:Sapprox}
\end{align}
Also, since $\int_{[-\delta,\delta]^d} g_\delta({\bf s})d{\bf s} = 1$, we may rewrite the definition of $S$ in (\ref{eq:Sdef}) as
\begin{align*}
S &=\sum_{\substack{{\bf k}\in\Z^d\\ K{\bf k}\in B_T}} \int_{[-\delta, \delta]^d} g_\delta({\bf s}) f(x_0 u(K{\bf k}))d{\bf s}
\end{align*}
and combining this with (\ref{eq:Sapprox}), we obtain
\begin{align*}
|S_{{\rm  approx}}-S| &\leq \left(\sum_{\substack{{\bf k}\in\Z^d\\ K{\bf k}\in B_T}} \int_{[-\delta, \delta]^d} g_\delta({\bf s}) |f(x_0 u(K{\bf k}+{\bf s}))-f(x_0 u(K{\bf k}))|d{\bf s}\right) + |r(T,K,f,d)|.
\end{align*}
But note that from property (\ref{Sob4}) of Sobolev norms, we have 
\[
|f(x_0 u(K{\bf k}+{\bf s}))-f(x_0 u(K{\bf k}))|\ll \sob{\infty}{1}(f) |{\bf s}| \ll \sob{\infty}{1}(f) \delta
\]
 for ${\bf s} \in [-\delta, \delta]^d.$
Together with our error bound, this implies that
\begin{align*}
|S_{{\rm  approx}} - S| &\ll \left( \sum_{\substack{{\bf k}\in\Z^d\\ K{\bf k}\in B_T}}\int_{[-\delta,\delta]^d} g_\delta({\bf s})d{\bf s}\right)\delta\sob{\infty}{1}(f)  + T^{d-1}\frac{\max_i K_i}{|K|}\sob{\infty}{0}(f)\\
&= \left(\#\{{\bf k}\in\Z^d | K{\bf k}\in B_T\} \delta +T^{d-1}\frac{\max_i K_i}{|K|}\right)\sob{\infty}{1}(f)
\end{align*}
once again, because $\int_{[-\delta,\delta]^d} g_\delta({\bf s})d{\bf s} = 1$.
But $\#\{{\bf k}\in\Z^d | K{\bf k}\in B_T\} \approx  T^d/|K|$, also with an error of magnitude $\ll T^{d-1}\max_i K_i/|K|$ (for reasons analagous to those illustrated in Figure \ref{fig:suminterror}).  Therefore
\begin{align}
|S_{{\rm  approx}} - S| \ll \left(\frac{T^d}{|K|}\delta+\frac{T^{d-1}\max_i K_i}{|K|}\right)\sob{\infty}{1}(f).\label{eq:SSapprox}
\end{align}

To show that $S_{{\rm  approx}}$ and $\dst{\frac{T^d}{|K|}\int_X f dm_X}$ are close, we observe that by Poisson summation,
\begin{align}
\sum_{{\bf k}\in\Z^d}  g_\delta({\bf t} - K{\bf k}) &= \sum_{{\bf k}\in\Z^d}  g_\delta({\bf t} + K{\bf k})\nonumber\\
&= \sum_{{\bf k}\in\Z^d} \widetilde{g_\delta}( K^{-1}{\bf t} + {\bf k})\nonumber\\
&= \sum_{{\bf k}\in\Z^d} \psi_{K^{-1}\bf k}({\bf t}) \widehat{\widetilde{g_\delta}}({\bf k})\label{eq:expakSum}
\end{align}
where $\psi_{K^{-1}\bf k}({\bf t}) = e^{2\pi i {\bf k} \cdot (K^{-1}{\bf t})} = e^{2\pi i (K^{-1}{\bf k}) \cdot {\bf t}}$ and $\widehat{\widetilde{g_\delta}}$ is the multivariate Fourier transform of $\widetilde g_\delta({\bf x}) = g_\delta(K{\bf x})$.  When we substitute (\ref{eq:expakSum}) into the definition of $S_{\rm  approx}$ given in (\ref{eq:Sapproxdef}), we get
\begin{align*}
S_{{\rm  approx}} &= \int_{B_T} \left(\sum_{{\bf k}\in\Z^d} \psi_{K^{-1}\bf k}({\bf t}) \widehat{\widetilde{g_\delta}}({\bf k})\right) f(x_0 u({\bf t}))d{\bf t}\\
&=\sum_{{\bf k}\in\Z^d} \widehat{\widetilde{g_\delta}}({\bf k})\left(\int_{B_T} \psi_{K^{-1}\bf k}({\bf t}) f(x_0 u({\bf t}))d{\bf t}\right)
\end{align*}
where Fubini's Theorem allows us to switch the order of the sum and the integral.  Similarly,
\begin{align*}
\frac{T^d}{|K|}\int_X f dm_X &= \left(\int_{B_T}\sum_{{\bf k}\in\Z^d}  g_\delta({\bf t} - K{\bf k}) d{\bf t}+\mathcal{O}\left(\frac{T^{d-1}\max_i K_i}{|K|}\right)\right)\int_X f dm_X\\
&=\sum_{{\bf k}\in\Z^d} \widehat{\widetilde{g_\delta}}({\bf k})\left(\int_{B_T} \psi_{K^{-1}\bf k}({\bf t}) \int_X f dm_X d{\bf t}\right) + \mathcal{O}\left(\frac{T^{d-1}\max_i K_i}{|K|}\sob{\infty}{0}(f)\right)
\end{align*}
where we have used that $|\int_X f dm_X|\leq \sob{\infty}{0}(f)$.  Thus
\begin{align}
\left|S_{{\rm  approx}}-\frac{T^d}{|K|}\int_X f dm_X\right| &= \left|\sum_{{\bf k}\in\Z^d} \widehat{\widetilde{g_\delta}}({\bf k})\int_{B_T} e^{2\pi i {\bf k} \cdot (K^{-1}{\bf t})}\left(f(x_0 u({\bf t})) - \int_X f dm_X\right) d{\bf t}\right|\nonumber\\
&\hspace{15pt}+\mathcal{O}\left(\frac{T^{d-1}\max_i K_i}{|K|}\sob{\infty}{0}(f)\right)\nonumber\\
&= \left|\sum_{{\bf k}\in\Z^d} \widehat{\widetilde{g_\delta}}({\bf k})|B_T|\mu_{T,\psi_{K^{-1}{\bf k}}}(f)\right|+\mathcal{O}\left(\frac{T^{d-1}\max_i K_i}{|K|}\sob{\infty}{0}(f)\right).\label{eq:SapproxIntegral}
\end{align}
Then since $R>C$, we can apply Lemma \ref{lem:VenkLem} to obtain
\begin{align*}
\left|S_{{\rm  approx}}-\frac{T^d}{|K|}\int_X f dm_X\right| &\ll_f T^d R^{-b}\sob{\infty}{\ell}(f) \sum_{{\bf k}\in\Z^d} \widehat{\widetilde{g_\delta}}({\bf k}) + \frac{T^{d-1}\max_i K_i}{|K|}\sob{\infty}{0}(f)
\end{align*}
(by direct computation we can see that $\widehat{\widetilde{g_\delta}}$ is positive).  Observe how it was crucial here that the result in Lemma \ref{lem:VenkLem} was uniform over characters.

Finally, again by Poisson summation, we have
\begin{align}
\sum_{{\bf k}\in\Z^d} \widehat{\widetilde{g_\delta}}({\bf k})
&= \sum_{{\bf k}\in\Z^d} \widetilde{g_\delta}({\bf k})\nonumber\\
&= \sum_{{\bf k}\in\Z^d} g_\delta(K{\bf k})\nonumber\\
&= g_\delta (0,\dots, 0)= \delta^{-d}\nonumber
\end{align}
since $\text{supp }(g_\delta) \subseteq [-\delta, \delta]^d$ and $\delta<\min_i K_i/2$ implies that $g_\delta(K{\bf k})= 0$ for ${\bf k} \neq (0, \dots, 0)$.  Substituting this into equation (\ref{eq:SapproxIntegral}), combining it with (\ref{eq:SSapprox}), and using property (\ref{Sob1}) of Sobolev norms, we get 
\begin{align*}
\left|S-\frac{T^d}{|K|}\int_X f dm_X\right| \ll \left(T^d R^{-b}\delta^{-d}+\frac{T^d}{|K|}\delta +\frac{T^{d-1}\max_i K_i}{|K|}\right)\sob{\infty}{\ell}(f).
\end{align*}
We can optimize the first two terms by choosing $\delta = (|K|/R^b)^{1/(d+1)}$.  Observe that our only resrtiction on $\delta$ was that $\delta < \min_i K_i/2$.  This will be achieved with our choice of $\delta$ so long as $R>(2/\min K_i)^{(d+1)/b}|K|^{1/b}$.
Thus, under these conditions,
\begin{align*}
\left|S-\frac{T^d}{|K|}\int_X f dm_X\right| \ll \left(T^d R^{-b/(d+1)} |K|^{-d/(d+1)} +\frac{T^{d-1}\max_i K_i}{|K|}\right)\sob{\infty}{\ell}(f).
\end{align*}

\end{proof}

If $K$ has all diagonal entries of equal weight (in abuse of notation, say all of weight $K$) then we get the following corollary which will be of use to us in the next section.

\begin{cor}
Let $T\geq K\geq 1$.  There exists a constant $\tilde C>0$ (depending only on $n$ and  $d$) such that for all $x_0\in X$ satisfying (\ref{eq:dio1}) with $T>R>\tilde C$, we have
\begin{align*}
\left|\sum_{\substack{{\bf k}\in\Z^d\\ K{\bf k}\in B_T}} f(x_0 u(K{\bf k}))-\frac{T^d}{K^d}\int_X f dm_X\right| \ll T^d R^{-b/(d+1)}K^{-d^2/(d+1)} \sob{\infty}{\ell}(f).
\end{align*}\label{cor:arithequidistcor}
\end{cor}

\begin{proof}
This is a straightforward application of the previous theorem, observing that in this case $(2/\min K_i)^{(d+1)/b}|K|^{1/b}=(2^{d+1}/K)^{1/b}\leq 2^{(d+1)/b}$ since $K\geq 1$.  Thus the theorem holds with $\tilde C =\max(C,2^{(d+1)/b})$.  Moreover, the second error term in Theoerem \ref{thm:arithequidist} in this case is simply $T^{d-1}K^{1-d}$, and since $K, R<T$ and $b<1$, this term decays more quickly than the first and can be ignored.
\end{proof}

\begin{rem}
For $X=\Gamma\backslash G$ where $\Gamma$ is a cocompact lattice, we have the following basepoint-independent versions of Lemma \ref{lem:VenkLem}, Theorem \ref{thm:arithequidist}, and Corollary \ref{cor:arithequidistcor}.
\begin{lem}
There exists $b>0$ (depending on $n$, $d$, and $\Gamma$) such that for all $T$ large enough, we have
\begin{align*}
\left|\mu_{T,\psi}(f)\right|\ll_\Gamma T^{- b} \sob{\infty}{\ell}(f)
\end{align*}
for any $f\in C^\infty(X)$, $x_0\in X$, and additive character $\psi$.\label{lem:VenkLemcmpt}
\end{lem}
\begin{thm}
Let $K = {\rm  diag}(K_1,\cdots, K_d)$ with $T\geq K_1, \dots, K_d\geq1$ and determinant $|K|$.  Then for all $T$ large enough, we have
\begin{align*}
\left|\sum_{\substack{{\bf k}\in\Z^d\\ K{\bf k}\in B_T}} f(x_0 u(K{\bf k}))-\frac{T^d}{|K|}\int_X f dm_X\right| \ll_\Gamma \left(T^{d-b/(d+1)} |K|^{-d/(d+1)} +\frac{T^{d-1}\max_i K_i}{|K|}\right)\sob{\infty}{\ell}(f)
\end{align*}
for all $f\in C^\infty(X)$ and $x_0\in X$.\label{thm:arithequidistcmpt}
\end{thm}
\begin{cor}
Let $T\geq K\geq 1$.  Then for all $T$ large enough, we have
\begin{align*}
\left|\sum_{\substack{{\bf k}\in\Z^d\\ K{\bf k}\in B_T}} f(x_0 u(K{\bf k}))-\frac{T^d}{K^d}\int_X f dm_X\right| \ll_\Gamma T^{d-b/(d+1)}K^{-d^2/(d+1)} \sob{\infty}{\ell}(f)
\end{align*}
for all $f\in C^\infty(X)$ and $x_0\in X$.\label{cor:arithequidistcorcmpt}
\end{cor}
The proofs of these results are completely analagous to the correspoding proofs for $\slnz\backslash\slnr$, but use the basepoint-independent equidistribution result stated in (\ref{eq:equidistcmpt}) instead of Theorem \ref{thm:equidist}.  As before, we may remove dependence on the lattice $\Gamma$ for $n\geq 3$ and for $n=2$ if $\Gamma$ is a congruence lattice.
\end{rem}


\section{Sieving and Orbits Along Almost-Primes}\label{sect:Sieving}

\subsection{$\Gamma$ Cocompact}

Let $\Gamma$ be a cocompact lattice in $G=\slnr$ and let $u({\bf t})$ be an abelian horospherical flow on $X=\Gamma/G$, as in Section \ref{sect:EquiArith}.  We know that the orbit of $u({\bf t})$ equidistributes with a uniform rate for all $x_0\in X$, 
and that as a consequence we have a uniform rate of equidistribution along multivariable arithmetic sequences of the form given in Corollary \ref{cor:arithequidistcorcmpt}.
Here and throughout this section, assume ${\bf k}=(k_1,\cdots, k_d)\in\Z^d$.

We want to understand the behavior of orbits at almost-prime entries of $u({\bf t})$.  More precisely, we want to understand averages of positive $f\in C_c^\infty(X)$ over points in $B_T$ that have entries with fewer than a certain fixed number of primes in their prime factorization.  

To investigate this question, we will use the combinatorial sieve from Theorem \ref{thm:CombSieve}.
In the context of our problem, we want to define
\begin{align*}
S(A,P) := \sum_{\substack{{\bf k}\in B_T\\ \gcd(k_1\cdots k_d, P)=1}} f(x_0u({\bf k}))
\end{align*}
where $f\in C^\infty_c(X)$, $f\geq 0$, and $P$ is the product of primes less than $z$ (to be determined). That is, we are summing over integer points in $B_T$ with entries containing no primes smaller than $z$. Then let
\begin{align*}
A = \{a_n\} := \left\{ \sum_{\substack{{\bf k}\in B_T\\ k_1\cdots k_d = n}} f(x_0u({\bf k})) \right\}
\end{align*}
and observe that
\begin{align*}
S_K(A,P) := \sum_{n\equiv 0 \mod K} \sum_{\substack{{\bf k}\in B_T\\ k_1\cdots k_d = n}} f(x_0u({\bf k}))
= \sum_{\substack{{\bf k}\in \tilde B_T\\ K|k_1k_2\cdots k_d}} f(x_0 u({\bf k}))
\end{align*}
where $\tilde B_T = (0,T]^d$ (since the index $n$ starts at 1 we want to avoid counting terms of the form $K$ divides $0$).

\begin{figure}
\begin{center}
\begin{tikzpicture}[scale=.8]
\draw[fill=gray!20]  (0,0) -- (0,3) -- (3,3) -- (3,0) -- cycle;
\node [gray] at (.75,.5) {$B_K$};
\node [below] at (3,0) {$K$};
\node [below] at (6,0) {$2K$};
\node [below] at (9,0) {$3K$};
\node [below] at (10,0) {$T$};
\node [left] at (0,3) {$K$};
\node [left] at (0,6) {$2K$};
\node [left] at (0,9) {$3K$};
\node [left] at (0,10) {$T$};
\node [below] at (2,1.5) {${\bf \tilde k}$};
\draw [dashed, red] [-] (2,1.5) -- (12,1.5);
\draw [dashed, red] [-] (2,1.5) -- (2,11.5);
\draw [dashed, red] [-] (12,1.5) -- (12,11.5);
\draw [dashed, red] [-] (2,11.5) -- (12,11.5);
\draw [dashed, red] [-] (5,1.5) -- (5,11.5);
\draw [dashed, red] [-] (8,1.5) -- (8,11.5);
\draw [dashed, red] [-] (11,1.5) -- (11,11.5);
\draw [dashed, red] [-] (2,4.5) -- (12,4.5);
\draw [dashed, red] [-] (2,7.5) -- (12,7.5);
\draw [dashed, red] [-] (2,10.5) -- (12,10.5);
\draw [-] (0,0) -- (10,0);
\draw [-] (0,0) -- (0,10);
\draw [-] (10,0) -- (10,10);
\draw [-] (0,10) -- (10,10);
\draw [-] (3,0) -- (3,10);
\draw [-] (6,0) -- (6,10);
\draw [-] (9,0) -- (9,10);
\draw [-] (0,3) -- (10,3);
\draw [-] (0,6) -- (10,6);
\draw [-] (0,9) -- (10,9);
\foreach \Point in {(.5,3),(.5,6),(.5,9),
(1,1.5),(1,3),(1,4.5),(1,6),(1,7.5),(1,9),
(1.5,1),(1.5,2),(1.5,3),(1.5,4),(1.5,5),(1.5,6),(1.5,7),(1.5,8),(1.5,9),(1.5,10),
(2,1.5),(2,3),(2,4.5),(2,6),(2,7.5),(2,9),
(2.5,3),(2.5,6),(2.5,9),
(3,.5),(3,1),(3,1.5),(3,2),(3,2.5),(3,3),(3,3.5),(3,4),(3,4.5),(3,5),(3,5.5),(3,6),(3,6.5),(3,7),(3,7.5),(3,8),(3,8.5),(3,9),(3,9.5),(3,10),
(3.5,3),(3.5,6),(3.5,9),
(4,1.5),(4,3),(4,4.5),(4,6),(4,7.5),(4,9),
(4.5,1),(4.5,2),(4.5,3),(4.5,4),(4.5,5),(4.5,6),(4.5,7),(4.5,8),(4.5,9),(4.5,10),
(5,1.5),(5,3),(5,4.5),(5,6),(5,7.5),(5,9),
(5.5,3),(5.5,6),(5.5,9),
(6,.5),(6,1),(6,1.5),(6,2),(6,2.5),(6,3),(6,3.5),(6,4),(6,4.5),(6,5),(6,5.5),(6,6),(6,6.5),(6,7),(6,7.5),
(6,8),(6,8.5),(6,9),(6,9.5),(6,10),
(6.5,3),(6.5,6),(6.5,9),
(7,1.5),(7,3),(7,4.5),(7,6),(7,7.5),(7,9),
(7.5,1),(7.5,2),(7.5,3),(7.5,4),(7.5,5),(7.5,6),(7.5,7),(7.5,8),(7.5,9),(7.5,10),
(8,1.5),(8,3),(8,4.5),(8,6),(8,7.5),(8,9),
(8.5,3),(8.5,6),(8.5,9),
(9,.5),(9,1),(9,1.5),(9,2),(9,2.5),(9,3),(9,3.5),(9,4),(9,4.5),(9,5),(9,5.5),(9,6),(9,6.5),(9,7),(9,7.5),(9,8),(9,8.5),(9,9),(9,9.5),(9,10),
(9.5,3),(9.5,6),(9.5,9),
(10,1.5),(10,3),(10,4.5),(10,6),(10,7.5),(10,9)}
{
    \fill \Point circle[radius=2pt];
}
\foreach \Point in {(2,10.5),(5,10.5),(8,10.5),(11,10.5),(11,7.5),(11,4.5),(11,1.5)}
{
    \fill[red] \Point circle[radius=2pt];
}
\end{tikzpicture}
\caption{In $S_K(A,P)$ we are summing over the integer points in $\tilde B_T$ such that $K|k_1\cdots k_2$ (shown in black).  We may do this by summing over shifted grids based at each of the points in the first box $\tilde B_K$ (shaded in gray).  However, this introduces an error determined by $\sob{\infty}{0}(f)$ and the number of points in each of these shifted grids falling outside $B_T$ (shown in red).  The number of such points can be bounded by $T^{d-1}K^{1-d}$, as we have seen before.}\label{fig:firstboxsumerror}
\end{center}
\end{figure}
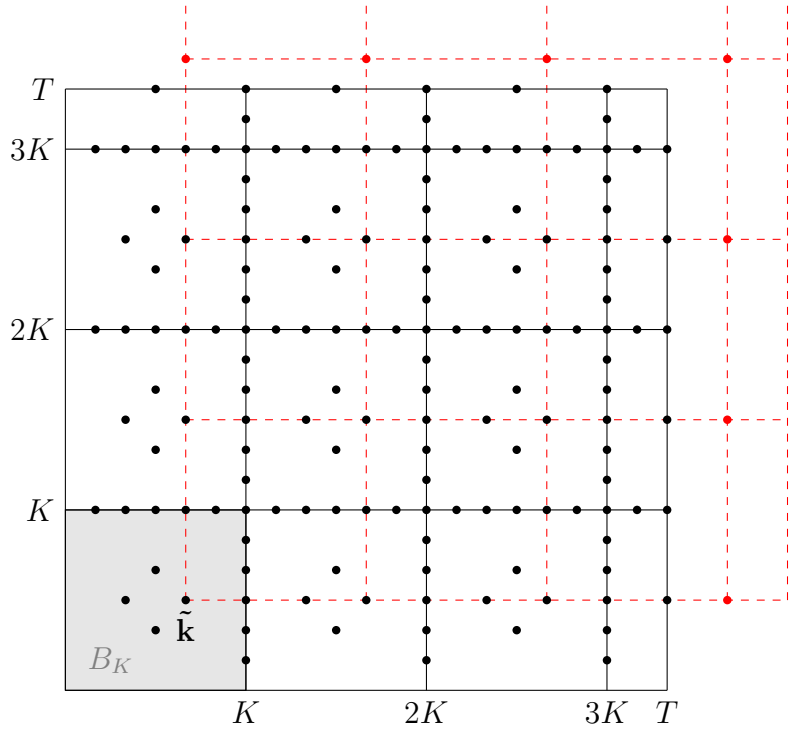

Notice that $K | k_1\cdots k_d$ if and only if $K | k_1\cdots (k_i + K) \cdots k_d$, that is, the collection of points that we are summing over is periodic with period $K$ in each coordinate.  Thus we can rewrite $S_K(A,P)$ as a sum over cubic grids of side length $K$ based at each point in the first box $B_K$:
\begin{align}
\sum_{\substack{{\bf k}\in \tilde B_T\\ K|k_1\cdots k_d}} f(x_0 u({\bf k})) = \sum_{\substack{{\bf \tilde k}\in \tilde B_K\\ K|\tilde k_1\cdots\tilde k_d}} \left( \sum_{K{\bf k} \in B_T} f(x_0 u({\bf \tilde k}) u(K{\bf k})) + \mathcal{O}(T^{d-1}K^{1-d}\sob{\infty}{0}(f))\right)\label{eq:doubleSum}
\end{align}
where the error arises from the fact that a point ${\bf \tilde k}+K{\bf k}$ for ${\bf k} \in B_T$ may, in fact, fall outside of $B_T$ (see Figure \ref{fig:firstboxsumerror}).

From Corollary \ref{cor:arithequidistcorcmpt}, we know that at each basepoint $x_{\bf \tilde k} = x_0 u({\bf \tilde k})$, we have
\begin{align}
\sum_{K{\bf k}\in B_T} f(x_{\bf \tilde k} u(K{\bf k})) = \frac{T^d}{K^d} \int f dm_X +\mathcal{O}\left(T^{d-b/(d+1)}K^{-d^2/(d+1)}\sob{\infty}{\ell}(f)\right).\label{eq:whateva}
\end{align}
If we let
\begin{align}
G_d (K) := \#\{ {\bf  k} \in \tilde B_K |  k_1 \cdots  k_d \equiv 0 \mod K \}\label{eq:Gddef}
\end{align}
then (\ref{eq:doubleSum}) together with (\ref{eq:whateva}) says that
\begin{align*}
S_K(A,P) = \sum_{\substack{{\bf k}\in \tilde B_T\\ K|k_1k_2\cdots k_d}} f(x_0 u({\bf k})) = \frac{G_d(K)}{K^d} \mathcal{X} +  r(f,K,T)
\end{align*}
where $\mathcal{X} = T^d \int f dm_X$ and
\begin{align*}
| r(f,K,T)|&\ll G_d(K)T^{d-b/(d+1)}K^{-d^2/(d+1)}\sob{\infty}{\ell}(f).
\end{align*}
This suggests that our function $g(K)$ in Theorem \ref{thm:CombSieve} should be $G_d(K)/K^d$, but it remains to show that this function satisfies the sieve axioms (\ref{axiom1}) and (\ref{axiom3}) and that the corresonding remainders satisfy the condition in axiom (\ref{axiom2}) for appropriately chosen $D(\mathcal{X})=D(T)$.  We will start with a lemma outlining some of the properties of the function $G_d$, the proof of which is given in Appendix \ref{app:Gd}.

\begin{lem}\label{lem:Gd}
 For any integers $K,d\geq 1$, the following hold:
\begin{enumerate}[(i)]
\item (Iterated sum formula) 
\[
G_d(K) = \sum^K_{k_{d-1} = 1}\cdots \sum^K_{k_1 = 1} \gcd(K, k_1\cdots k_{d-1}).
\]\label{Gdprops1}
\item (Recursive formula) Let ${\rm  Id}^d(K) = K^d$.  Then
\[
G_{d+1} = {\rm  Id}^d\ast(\phi\cdot G_d).
\]\label{Gdprops2}
\item $G_d$ is multiplicative.\label{Gdprops3}
\item (Behavior at primes) Let $p$ be a prime.  Then 
\[
G_d(p) = p^d-(p-1)^d.
\]\label{Gdprops4}
\item (Dirichlet series bound) For real $x>e$ and $s<d$, 
\[
\sum_{K\leq x} \frac{G_d(K)}{K^s} \ll_{s,d} x^{d-s}(\log x)^{d-1}.
\]\label{Gdprops5}
\end{enumerate}
\end{lem}

\begin{rem}
By convention, we let $G_1(K) = \gcd(K,1) = 1 = \#\{0<k\leq K | k \equiv 0 \mod K\}$, which trivially satisfies the above properties, so long as the empty product in (\ref{Gdprops1}) is properly interpreted to be 1.
\end{rem}

\begin{rem}
Notice that $G_2(K) = \sum^K_{j=1} \gcd(K, j)$ is Pillai's arithmetical function,\footnote{ The values of $G_2(K)$ for $K=1, 2, 3, \dots$ are given as sequence A018804 in the OEIS (see \cite{OEIS}).}
 a multiplicative function first considered by Ces\`{a}ro and rediscovered by Pillai in \cite{Pillai}.  For this function, property (\ref{Gdprops2}) is the well-known identity $G_2 = {\rm  Id}\ast\phi$.\footnote{ For $G_2$, much else is known.  In terms of Dirichlet convolution, we also have the useful identity $G_2 = \mu \ast ({\rm  Id}\cdot \tau)$, where $\mu$ is the M\"{o}bius function and $\tau$ is the divisor function.  In \cite{Broug1}, Broughan used this to derive a closed form for the Dirichlet series in terms of the Riemann zeta function, as well as an asymptotic formula for the partial sums of the Dirichlet series.  The asymptotics for partial sums of the Dirichlet series were later refined by \cite{BordPillai}, \cite{Broug2}, and \cite{Tanigawa}.}
As noted in \cite{OEIS}, $G_2(K)$ counts the number of non-congruent solutions to the equation $k_1 k_2 \equiv 0 \mod K$.  From the definition of $G_d$ in (\ref{eq:Gddef}), we can similarly see that $G_d(K)$ counts the number of non-congruent solutions to $k_1k_2\cdots k_d \equiv 0 \mod K$, so in this way $G_d$ can be considered a generalization of Pillai's arithmetical function.\footnote{ Other generalizations of Pillai's arithmetical function have been studied.  Examples include \cite{LongGcdGeneral}, \cite{TothGcdGeneral}, \cite{BordGcdGeneral}, \cite{HaukGcdGeneral}, and \cite{Toth}, however none of these include the generalization given here.  In \cite{TothSimilar}, T\`{o}th considers a generalization that is very similar to ours, and in the notation of that paper, $G_d(K) = A_{d-1}(K)K^{d-1}$.  Lemma \ref{lem:Gd} (\ref{Gdprops2}) and (\ref{Gdprops3}) can thus be considered corollaries of results proved in \cite{TothSimilar}, but we prove them in Appendix \ref{app:Gd} in order to keep the paper self-contained.  T\`{o}th also gives a formula for the Dirichlet series of this generalization in terms of the Dirichlet series of a related arithmetic function, however we will need an explicit estimate for the partial sums of the Dirichlet series where it does not converge, which we develop as property (\ref{Gdprops5}) of Lemma \ref{lem:Gd}.}
In addition to those presented here, there are undoubtedly many other useful properties and interpretations of the generalized Pillai's functions $G_d$, which could be an interesting area of future study.\footnote{ For example, for $K$ squarefree, we have the bound $G_d(K) \leq K^{d-1} d^{\omega(K)}$, where $\omega(K)$ counts the number of (distinct) primes dividing $K$ (since $K$ is squarefree, we have $\omega(K)=\Omega(K)$).  This bound can be derived from the formula for primes along with multiplicativity and can be used along with known estimates for $\omega(K)$ as an alternative to Lemma \ref{lem:Gd} (\ref{Gdprops5}) to verify sieve axiom (\ref{axiom2}).}

\end{rem}

We can now verify that the sieve axioms in Theorem \ref{thm:CombSieve} are satisfied, which gives us the following theorem.

\begin{thm}
Let $u$ be a $d$-dimensional abelian horospherical flow on $X=\Gamma\backslash\slnr$ for $\Gamma$ cocompact, and let $P$ be the product of primes less than $T^\alpha$ for $\alpha<b/ 9d^2$, where $b$ is the constant from Lemma \ref{lem:VenkLemcmpt}.  Then for any $x_0 \in X$, positive $f\in C^\infty(X)$, and $T$ large enough (depending on $n$, $d$, $\Gamma$, and $f$), we have
\[
\sum_{\substack{{\bf k}\in B_T\\ \gcd(k_1\cdots k_d, P) = 1}} f(x_0u({\bf k})) \asymp_{\Gamma} \left(\frac{T}{\log T}\right)^d\int f dm_X.
\]
\label{thm:SAPcmpt}
\end{thm}

\begin{rem}
It is not clear from the statement of Theorem \ref{thm:CombSieve}, but from \cite{CombSieve} it can be found that the dependence on $f$ arising from the implicit constant in sieve axiom (\ref{axiom2}) can be entirely absorbed by the implicit constant determining how large we require $T$ to be to get the result.  As usual, the implict constant in the conclusion of this theorem depends also on $n$ and $d$, and dependence on $\Gamma$ may be removed if $n\geq 3$ or if $\Gamma$ is a congruence lattice.
\end{rem}

\begin{rem}
Let $\phi(x,y)$ be the number of positive integers $\leq x$ not divisible by any prime $\leq y$ for $x\geq y\geq 2$.  It is known that
\[
\phi(x,y) = \frac{x\omega(\log x/\log y)-y}{\log y} + \mathcal{O}\left(\frac{x}{(\log y)^2}\right)
\]
where $\omega:[1,\infty)\to[1/2,1]$ is the Buchstab function.  
Thus, the number of integers in $[0,T]$ not divisible by any prime less than $T^\alpha$ for $\alpha<1$ is given by
\[
\phi(T,T^\alpha) = \frac{\omega(1/\alpha)T}{\alpha\log T} - \frac{T^\alpha}{\alpha\log T} + \mathcal{O}\left(\frac{T}{(\alpha \log T)^2}\right).
\]
Thus the number of points ${\bf k}\in B_T$ such that $\gcd(k_1\cdots k_d, P) = 1$ where $P$ is the product of primes less than $T^\alpha$ is $\phi(T,T^\alpha)^d$, which grows asymptotically like $(T/\log T)^d$ as $T\to \infty$.  Although our result above only states that there is an upper and lower bound with respect to this quantity, it hints that there may be underlying equidistribution behavior.
\end{rem}

\begin{proof}

We need to show that sieve axioms (\ref{axiom1}), (\ref{axiom2}), and (\ref{axiom3}) are satisfied for 
\[
S_K(A,P) := \sum_{\substack{{\bf k}\in \tilde B_T\\K|k_1\cdots k_d}} f(x_0 u({\bf k})) = g(K) \mathcal{X} + r(f,K,T)
\]
where $g(K) = G_d(K)/K^d$, $\mathcal{X}=T^d\int f dm_X$, and 
\[| r(f,K,T)|\ll G_d(K)T^{d-b/(d+1)}K^{-d^2/(d+1)} \sob{\infty}{\ell}(f)
\]
 (see discussion at the beginning of the section).

To verify sieve axiom (\ref{axiom1}), note that since $G_d(K)$ is multiplicative by Lemma \ref{lem:Gd} (\ref{Gdprops3}), it follows that $g(K) = G_d(K)/K^d$ is multiplicative.  Furthermore, by Lemma \ref{lem:Gd} (\ref{Gdprops4}), we know that at primes
\begin{align*}
0< g(p) &= \frac{p^d-(p-1)^d}{p^d}\\
&= 1 - \left(\frac{p-1}{p}\right)^d\\
&\leq 1 - \left(\frac{2-1}{2}\right)^d = 1-\frac{1}{2^d}
\end{align*}
since $p\geq2$.  So sieve axiom (\ref{axiom1}) is satisfied with, for example, $c_1 = 2^{d+1}$.

To verify sieve axiom (\ref{axiom2}), observe that by Lemma \ref{lem:Gd} (\ref{Gdprops5}), we know that for $D$ large enough,
\begin{align*}
\sum_{K<D} |\tilde r(f,K,T)| &\ll T^{d-b/(d+1)}\sob{\infty}{\ell}(f)\sum_{K<D}\frac{G_d(K)}{K^{d^2/(d+1)}}\\
&\ll_f   T^{d-b/(d+1)} D^{d/(d+1)}( \log D)^{d-1}.
\end{align*}
Now if we let $D=T^\eta$ for any $\eta<b/d$, say $\eta=b/d - 2(d+1)\epsilon$ for some $\epsilon>0$, then we have
\begin{align*}
\sum_{K<D} |\tilde r(f,K,T)| &\ll_{f,\epsilon}  T^{d-2d\epsilon} ( \log T)^{d-1}.
\end{align*}
But since $\log T$ asymptotically grows more slowly than any positive power of $T$, we can say that for $T$ large enough, $\log T \ll_{\epsilon} T^{d\epsilon/(d-1)}$.  Hence
\[
\sum_{K<D} |\tilde r(f,K,T)| \ll_{f,\epsilon} T^{d-2d\epsilon}(T^{d\epsilon/(d-1)})^{d-1} = (T^d)^{1-\epsilon} \ll_f \mathcal{X}^{1-\epsilon}.
\]

To verify sieve axiom (\ref{axiom3}), notice that (by the binomial theorem)
\begin{align}
g(p) = \frac{p^d-(p-1)^d}{p^d} = \frac{d}{p} - \sum_{i=2}^d \frac{a_i}{p^i} \label{eq: A31}
\end{align}
where $a_i = (-1)^i {d \choose i}$. Since $\sum_{j=1}^\infty\log(j)/j^i$ converges for any $i>1$, we have that 
\begin{align}
\left|\sum_{w\leq p \leq z} \sum_{i=2}^d \frac{a_i\log p}{p^i}\right| 
\leq \sum_{i=2}^d|a_i|\sum_{j=1}^\infty \frac{\log j}{j^i} = C_2 \label{eq: A32}
\end{align}
and by a corollary of the Prime Number Theorem, we know that
\begin{align*}
\sum_{p\leq x} \frac{\log p}{p} = \log(x) + \mathcal{O}(1).
\end{align*}
Hence,
\begin{align*}
\sum_{p\leq z} \frac{\log p}{p} - \sum_{p< w} \frac{\log p}{p} &= \log(z)-\log(w) + \mathcal{O}(1)\\
\sum_{w\leq p \leq z} \frac{\log p}{p}&= \log\frac{z}{w} +\mathcal{O}(1)
\end{align*}
i.e., there exists $C_2'$ such that
\begin{align}
\left| \sum_{w\leq p \leq z} \frac{\log p}{p} - \log\frac{z}{w}\right| \leq C_2' \label{eq: A33}
\end{align}
for all $2\leq w\leq z$.  Putting (\ref{eq: A31}), (\ref{eq: A32}), and (\ref{eq: A33})   together, we see that 
\begin{align*}
\left| \sum_{w\leq p \leq z} g(p) \log p - d\log\frac{z}{w}\right| &\leq d \left| \sum_{w\leq p \leq z} \frac{\log p}{p} - \log\frac{z}{w}\right| + \left|\sum_{w\leq p \leq z} \sum_{i=2}^d \frac{a_i\log p}{p^i}\right|\\
&\leq dC_2'+C_2
\end{align*}
which shows that axiom (\ref{axiom3}) is satisfied with sieve dimension $r=d$ and $c_2 = dC_2'+C_2$.

Since we have demonstrated that sieve axioms (\ref{axiom1}), (\ref{axiom2}), and (\ref{axiom3}) hold, we have the conclusion of Theorem \ref{thm:CombSieve}, which implies our result.
\end{proof}

Notice that if an integer $k<T$ has no prime factors less than $T^\alpha$, then it must have fewer than $1/\alpha$ prime factors total.  Hence, if we take $f$ to be a positive function supported on any small neighborhood, Theorem \ref{thm:SAPcmpt} tells us that we can take $T$ large enough so that averaging $f$ over integer points in $B_T$ with no prime factors less than $T^\alpha$ has a positive lower bound.  This means that the set of $(1/\alpha)$-almost-prime times hitting any neighborhood is nonempty, which gives us the theorem from the introduction with $M=1/\alpha$.

\begin{cor*}[Theorem \ref{thm:CmptThm}]
Let $u({\bf t})$ be an abelian horospherical flow of dimension $d$ on $X=\Gamma\backslash \slnr$ for $\Gamma$ cocompact.
Then there exists a constant $M$ (depending only on $n$, $d$, and $\Gamma$) such that for any
$x_0\in X$, the set 
\[
\{x_0u(k_1, k_2, \cdots, k_d) \hspace{2pt}\vert\hspace{2pt} k_i\in\Z \text{ has fewer than } M \text{ prime factors}\}
\]
is dense in $X$.
\end{cor*}


\subsection{The Space of Lattices}

Now consider $X = \Gamma\backslash G$ for the non-cocompact lattice $\Gamma=\slnz$. Since we no longer have a uniform rate of equidistribution for our abelian horospherical flow $u({\bf s})$, we will consider a basepoint $x_0= \Gamma g_0 \in X$ satisfying a Diophantine condition of the following form.

\begin{defn}
We say that $x=\Gamma g$ is \textit{strongly polynomially} $\delta$\textit{-Diophantine} if there exists a sequence $T_i\to\infty$ as $i\to\infty$ such that 
\[
\inf_{\substack{w\in \Lambda^j(\Z^n)\setminus\{0\}\\j=1, \cdots, n-1}} \sup_{{\bf t}\in [0,T_i]^d}\nn{wgu({\bf t})}>T_i^\delta
\]
for all $i\in\N$.\label{def:Dio}
\end{defn}

The motivation for this definition is that, as in the compact setting, we will want to apply sieving to learn about integer points haivng few prime factors. However, unlike in the compact case, we do not have a uniform rate of equidistribution, so we must consider the effect of the basepoint.  For a given time-scale $T$, to obtain information about almost-primes of a certain order, we would want $R$ in the basepoint condition (\ref{eq:dio1}) to look like a small power of $T$ (say $T^\delta$).  However, a theorem like that of Theorem \ref{thm:SAPcmpt} will require $T$ be ``large enough," which depends on the function $f$, and so any fixed time-scale $T$ is insufficient.  Moreover, the constant $\delta$ we are able to take at one time-scale may not work for a different time-scale, which affects the number of prime factors we allow for our almost-prime points.  The condition given in Definition \ref{def:Dio} ensures that for any function (hence any neighborhood in $X$) we will be able to find a time-scale large enough so that our sieving provides positive information about almost-primes of the same, fixed order.

Before moving on to the main theorem of this section, we briefly remark that this definition is a meaningful one.  In view of results in \cite{KMLogLaws}, we see that not only do such points exist, but any generic point for the flow $u$ will satisfy this definition for some positive $\delta$.

\begin{thm}\label{thm:SAPnoncmpt}
Let $u$ be an abelian horospherical flow on $X=\slnz\backslash\slnr$ and let $P$
be the product of primes less than $T^\alpha$ for $\alpha<\delta bn/9d(d^2+bn\kappa)$, where $b$ is the constant from Lemma \ref{lem:VenkLem} and $\kappa = \min(m,n-m)$ for $m$ as in (\ref{eq:U}).  Furthermore, let $x_0\in X$ be strongly polynomially $\delta$-Diophantine.  Then for any positive $f\in\CX$ there exists a sequence $T_i\to\infty$ as $i\to\infty$ where
\[
\sum_{\substack{{\bf k}\in B_{T_i}\\ \gcd(k_1\cdots k_d, P) = 1}} f(x_0u({\bf k})) \asymp
\left(\frac{T_i}{\log T_i}\right)^d\int f dm_X.
\]
\end{thm}

\begin{proof}
Let $f\in\CX$, $f\geq 0$, and let $u$ be an abelian horospherical flow as given in (\ref{eq:U}) of Section \ref{sect:EquiArith}.  
As in the compact setting, we want to use our equidistribution theorem for arithmetic sequences to say that 
\begin{align}
S_K(A,P):&=\sum_{\substack{{\bf k}\in B_{T_i}\\ K|k_1\cdots k_d}} f(x_0u({\bf k}))\nonumber\\
&= \sum_{\substack{{\bf \tilde k}\in \tilde B_K\\ K|\tilde k_1\cdots\tilde k_d}} \left( \sum_{K{\bf k} \in B_{T_i}} f(x_0u({\bf \tilde k})u(K{\bf k})) + \mathcal{O}(T_i^{d-1}K^{1-d})\right)\label{eq:doublesumnoncmpt}\\
&= g(K) \mathcal{X} + r(f,K,T_i)\nonumber
\end{align}
where $g(K)=G_d(K)/K^d$, $\mathcal{X}=T_i^d \int f dm_X$, and the error terms can be suitably controlled.
Unfortunately, we cannot apply the same equidistribution result to the shifted basepoints $x_0u({\bf \tilde k})$ since they will not necessarily satisfy the same Diophantine condition.  However, since $K$ is understood to be small in comparison to the $T_i$, all of the points in $B_K$ lie comparatively close to $x_0$.  Then since the Diophantine property varies continuously, we expect the points in this region to satisfy a Diophantine condition not much worse than that of $x_0$, and in fact we can make this quantitative.

Observe that if $x_0$ is strongly polynomially $\delta$-Diophantine, it means that condition (\ref{eq:dio1}) holds for the sequence of parameters $T=T_i$ and $R=T_i^{\delta/q}$, where $q=\sum_{\lambda_i<0}-m_i\lambda_i = d/n$ for abelian $u$ of this form.
That is, for all $j\in\{1,\cdots,n-1\}$ and $w\in \Lambda^j(\Z^n)\setminus\{0\}$, we have
\begin{align}
\exists\hspace{2pt} {\bf t}\in [0,T_i]^d \text{ s.t. } \nn{wg_0u({\bf t})}=\nn{wg_0u({\bf \tilde k})u({\bf t})u({\bf -\tilde k})}\geq T_i^\delta. \label{eq:conjcond}
\end{align}
Recall that any $w \in \Lambda^j(\R^n)$ can be written as a sum $w = \sum_I w_I e_I$ over multi-indices $I = (i_1, \cdots, i_j)$ with $0<i_j<\cdots<i_1<n$, coefficients $w_I\in\R$, and basis elements  $e_I = e_{i_1}\wedge\cdots\wedge e_{i_j}$ where $\{e_i\}_{1\leq i \leq n}$ is the standard basis on $\R^n$.  Recall also that the norm above is defined by
\begin{align*}
\nn{w} = \max_I |w_I|
\end{align*}
and that $G$ acts linearly on $\Lambda^j(\R^n)$ by sending a basis vector $e_{i_1}\wedge \cdots\wedge e_{i_j}$ to
\begin{align*}
(e_{i_1}\wedge \cdots\wedge e_{i_j}) g = (e_{i_1}g)\wedge \cdots\wedge(e_{i_j}g).
\end{align*}

Since our abelian horospherical subgroup has the form given in (\ref{eq:U}), we can write an arbitrary $u\in B_K^{-1} = u([-K,0]^d)$ as
\begin{align}
u =
\left(\begin{matrix}
&&\vline& a_{1(m+1)}&\cdots&a_{1n}\\
& I_{m}  & \vline& \vdots& & \vdots\\
&&\vline& a_{m(m+1)}&\cdots&a_{mn}\\
\hline
&0&\vline&&I_{n-m}&\\
\end{matrix} 
\right)
\end{align}
where $a_{ij}\in[-K,0]$ for all $1\leq i \leq m$ and $m+1\leq j \leq n$.  One may verify that 
\begin{align*}
e_i u = e_i + a_{i(m+1)} e_{m+1} + \cdots + a_{in} e_n
\end{align*}
for $1\leq i \leq m$, and 
\begin{align*}
e_i u = e_i
\end{align*}
for $m+1\leq i \leq n$.  Hence, when we take wedge products $(e_{i_1}u)\wedge \cdots\wedge(e_{i_j}u)$, we cannot get a coefficient of order greater than $K^m$, since only the first $m$ transformed basis vectors have nontrivial coefficients and none of these coefficients have magnitude greater than $K$.  On the other hand, we cannot get a coefficient of order larger than $K^{n-m}$, since only the basis vectors $e_{m+1}$ through $e_{n}$ carry nontrivial coefficients.  Thus if we let $\kappa:=\min\{m, n-m\}$, we find that
\begin{align*}
\nn{(e_{i_1}u)\wedge \cdots\wedge(e_{i_j}u)}\ll K^\kappa.
\end{align*}
Then for general $w\in \Lambda^j(\R^n)$ and $u\in B_K$, we have
\begin{align*}
\nn{wu}\ll K^\kappa\nn{w}.
\end{align*}
Thus from (\ref{eq:conjcond}), we can say that for any $w \in \Lambda^j(\Z^n)\setminus\{0\}$, $j\in\{1,\cdots, n-1\}$, there exists $ {\bf t}\in [0,T_i]^d$ such that 
\begin{align*}
K^\kappa \nn{wg_0u({\bf \tilde k})u({\bf t})} &\gg\nn{wg_0u({\bf \tilde k})u({\bf t})u({\bf -\tilde k})}\geq T_i^\delta
\end{align*}
so
\begin{align*}
\nn{wg_0u({\bf \tilde k})u({\bf s})} &\gg T_i^\delta/K^\kappa
\end{align*}
That is, for any $u({\bf \tilde k})\in B_K$, the shifted basepoint $x_0u({\bf \tilde k})$ satisfies a  Diophantine condition of the form (\ref{eq:dio1}) with new parameter proportional to $(T_i^\delta/K^\kappa)^{1/q}=(T_i^\delta/K^\kappa)^{n/d}$.  From Corollary \ref{cor:arithequidistcor}, this implies that for $T_i$ large enough (i.e. for $i$ large enough), we have equidistribution with \begin{align}
\sum_{K{\bf k}\in B_{T_i}} f(x_0 u(K{\bf k})) = \frac{T_i^d}{K^d}\int_X f dm_X+\mathcal{O}_f(T_i^d(T_i^\delta/K^\kappa)^{-nb/d(d+1)}K^{-d^2/(d+1)})\label{eq:arithEquiNoncmpt}
\end{align}
for any ${\bf k}\in B_K$.
Using this in (\ref{eq:doublesumnoncmpt}), we find that
\[
S_K(A,P) = g(K)\mathcal{X}+r(f,K,T_i)
\]
where
\[
|r(f,K,T_i)|\ll G_d(K)T_i^{d-\delta nb/d(d+1)}K^{(\kappa nb-d^3)/d(d+1)}\sob{\infty}{\ell}(f).
\]
Since we have already shown that the function $g(K) = G_d(K)/K^d$ satisfies sieve axioms (\ref{axiom1}) and (\ref{axiom3}) with sieve dimension $d$, it remains to verify sieve axiom (\ref{axiom2}).
From Lemma \ref{lem:Gd} (\ref{Gdprops5}), we know that
\begin{align*}
\sum_{K<D} |r(f,K,T_i)| &\ll_f T_i^{d-\delta nb/d(d+1)}\sum_{K<D} G_d(K)/K^{(d^3-\kappa nb)/d(d+1)}\\
&\ll T_i^{d-\delta nb/d(d+1)} D^{(d^2+\kappa nb)/d(d+1)}(\log D)^{d-1}.
\end{align*}
Now let $D = T^\eta$ for $\eta < \delta bn/(d^2+\kappa bn)$, say $\eta = \delta bn/(d^2+\kappa bn) - 2\epsilon d^2(d+1)/(d^2+\kappa nb)$.  As before, we know that for large enough $T_i$, $\log T_i \ll_\epsilon T_i^{d\epsilon/(d-1)}$.  This is enough to ensure that for $T_i$ large enough, the errors satisfy
\[
\sum_{K<D} |r(f,K,T_i)| \ll_{f,\epsilon} T^{d(1-\epsilon)}\ll_f \mathcal{X}^{1-\epsilon}.
\]
Thus for given $f$, the conclusion of Theorem \ref{thm:CombSieve} holds for all $i$ large enough, which gives us Theorem \ref{thm:SAPnoncmpt}.
\end{proof}

As before, if we consider positive $f$ supported on a neighborhood of $X$, the above theorem tells us that we may take $i$ large enough so that we have a positive lower bound on averages over almost-prime points with fewer than $1/\alpha = 9d(d^2+\kappa nb)/\delta bn$ prime factors, hence such points are dense in $X$.  This gives us the theorem for the space of lattices from the introduction.\footnote{ Strictly speaking, if we let $M_\delta = 1/\alpha$, this also has dependence on $\kappa$, which cannot be explicitly reduced to dependence on $n$ and $d$.  However, if we want to eliminate this dependence, we may replace $\kappa$ with $n/2$, since $\kappa\leq n/2$ in all cases.}

\begin{cor*}[Theorem \ref{thm:NoncmptThm}]
Let $u({\bf t})$ be an abelian horospherical flow of dimension $d$ on $X=\slnz\backslash\slnr$ and let $x_0\in X$ be strongly polynomially $\delta$-Diophantine for some $\delta>0$.  Then there exists a constant $M_\delta$ (depending on $\delta$, $n$, and $d$) such that 
\[
\{x_0u(k_1, k_2, \cdots, k_d) \hspace{2pt}\vert\hspace{2pt} k_i\in\Z \text{ has fewer than } M_\delta \text{ prime factors}\}
\]
is dense in $X$.
\end{cor*}


\section{Conclusion}\label{sect:concl}

In this paper we gave an effective equidistribution result for horospherical flows on the space of lattices and an effective rate of equidistribution for arithmetic sequences of entries in abelian horospherical flows on both the space of lattices and compact quotients of $\slnr$.  We then use sieve methods to derive an upper and lower bound for averages over almost-prime entries of abelian horospherical flows.  In the compact setting, we have as a result the density of integer entries having fewer than a fixed number of primes depending only on the dynamical system and not on the basepoint.  In the space of lattices, we consider the orbits of points satisfying a Diophantine condition with parameter $\delta$ and we prove the density of integer entries having fewer than a fixed number of primes depending on the system and on $\delta$.

There are several improvements and generalizations of this work that can be readily imagined. It seems likely that the methods used here can be generalized to quotients of 
connected, semisimple Lie groups by lattices.  It also seems possible that methods similar to those used in \cite{SarnackUbis} could be adapted to remove dependence on the basepoint in the noncompact case, yielding a uniform result for the density of almost-primes in the orbits of any generic point.
One could also generalize from abelian horospherical flows to arbitrary horospherical flows.  This could make the character analysis in Section \ref{sect:EquiArith} on arithmetic sequences more tricky, but it nonetheless seems doable.  Finally, the sieve methods used Section \ref{sect:Sieving} can be modified to learn about averages over points $(k_1,  \cdots, k_d)\in \R^d$  satisfying $\gcd(\mathcal{P}(k_1, \cdots, k_d), P) =1$, where $\mathcal{P}$ is a suitably nice irreducible polynomial (note that we considered the case where $\mathcal{P}(k_1, \cdots, k_d) = k_1k_2\cdots k_d$).


Of course, the more natural question is not what happens at almost-prime times, but what happens at prime times.  Unfortunately, it does not seem possible at present to use these methods to establish results about true primes, and additional ingredients or a wholly different approach may be required.  However, this result is significant in that it continues to lend support to the conjecture, already suggested by \cite{SarnackUbis}, that prime times in horospherical orbits are dense and possibly equidistributed.


\begin{appendix}

\section{Radius of Injection}\label{app:RadiusofInjection}

Let $G=\slnr$, $\Gamma=\slnz$, and $X=\Gamma\backslash G$ be the space of lattices.  We want to prove the following lemma for the radius of injection.

\begin{customlem}{\ref{lem:radiusofinjection}}
There exist constants $c_1, c_2>0$ (depending only on $n$) such that for any $0<\epsilon<c_1$, the projection map
\begin{align*}
\pi_x: B_r^G(e) &\to B_r^X(x)\\
g &\mapsto xg
\end{align*}
is injective for all $x\in L_\epsilon$, where $r=c_2\epsilon^n$.
\end{customlem}

To do this, we will first need some background on Siegel sets for the action of $\slnz$ on $\slnr$.

\subsection{Siegel Sets}\label{sect:SiegelSet}

Let $K = {\rm  SO}(n)$, let $A$ be the positive diagonal subgroup, and let $N$ be the subgroup of upper triangular unipotents.  The Iwasawa decomposition of $G$ is given by $G = NAK$.  One can use reduction theory for arithmetic groups to find a convenient way of writing $x\in X$ in terms of particular subsets of these subgroups.

Given $\epsilon>0$, define
\begin{align*}
A_\epsilon &= \left\{{\rm  diag}(a_1, \cdots, a_n)\in A \hspace{4pt} \vline \hspace{4pt} \frac{a_{i+1}}{a_{i}} \leq \epsilon \right\}\\
N_\epsilon &=\left\{ u \in N \hspace{2pt} \vline \hspace{4pt} |u_{i,j}| \leq \epsilon \hspace{6pt} \forall i<j \right\}.
\end{align*}
A Siegel set for $G$ is a set of the form $\Sigma_{s,t} := N_{s}A_{t} K$ for some $s,t>0$.

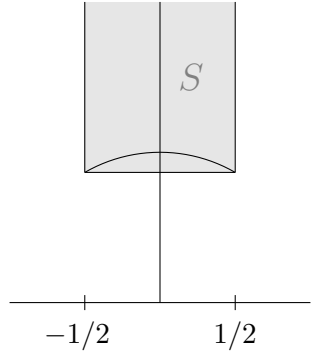
\begin{figure}
\begin{center}
\begin{tikzpicture}
\fill [gray!20] (1, 4) -- (1, 1.732) -- (-1, 1.732) -- (-1,4);
\node[right,gray] at (.1,3) {\large{$S$}};
\draw [-] (-2,0) -- (2,0);
\draw [-] (0,0) -- (0,4);
\draw [-] (1, 4) -- (1, 1.732);
\draw [-] (-1, 4) -- (-1, 1.732);
\draw [-] (-1, 1.732) -- (1, 1.732);
\draw [domain=60:120] plot ({2*cos(\x)}, {2*sin(\x)});
\node [below] at (-1.1,-.1) {\small{$-1/2$}};
\node [below] at (1,-.1) {\small{$1/2$}};
\draw [-] (-1,-.1)--(-1, .1);
\draw [-] (1, -.1)--(1, .1);
\end{tikzpicture}
\caption{The Siegel set $S = \Sigma_{1/2,2/\sqrt{3}}$ and the fundamental domain when $n=2$, represented in the Poincar\'e upper half plane.}
\label{fig:SiegelSet}
\end{center}
\end{figure}

Siegel sets can be thought of as a nice way of approximating a fundamental domain for the action of $\Gamma=\slnz$ on $G$
(see Figure \ref{fig:SiegelSet}).
This approximation can be optimized in the following sense:
For any 
$s\geq 1/2$ and $t\geq 2/\sqrt{3}$, $G=\slnr$ can be written as
\begin{align*}
G = \Gamma \Sigma_{s,t}
\end{align*}
(for details and a proof, see \cite{Rag} Theorem 10.4 or \cite{BekkaMayer} Theorem 5.1.7).

\subsection{Proof of Radius of Injection}

We start with the following well-known computation.

\begin{lem}\label{lem:Adg}
Let $g\in \Sigma_{\frac{1}{2},\frac{2}{\sqrt{3}}}$ satisfy $\Gamma g\in L_\epsilon$.  Then the operator norm of ${\rm  Ad}_g:{\rm  Mat}_{n\times n}(\R) \to {\rm  Mat}_{n\times n}(\R)$ satisfies
\begin{align*}
\nn{{\rm  Ad}_g}\ll \epsilon^{-n}
\end{align*}
where the implicit constant depends only on dimension $n$.
\end{lem}

\begin{proof}
Let $g=uak$ where $u\in U_{1/2}$, $a={\rm  diag}(a_1,\cdots, a_n)\in A_{2/\sqrt{3}}$, and $k\in {\rm  SO}(n)$. Let $\{e_i\}_{1\leq i \leq n}$ be the standard basis on $\R^n$ and fix $\nn{\cdot}$ to be the max matrix norm on ${\rm  Mat}_{n\times n}(\R)$ (any other norm will work equally well).

Notice that $e_n u = e_n$ for any $u\in U$ and that $e_n a = a_ne_n$ for $a\in A$.  Furthermore, since $k$ is an orthogonal matrix, we have that $\nn{v k} \leq \sqrt{n}\nn{v}$ for any $v\in \R^n$.  Then, since $\Gamma uak\in L_\epsilon$, we know that $\nn{vuak}\geq \epsilon$ for all $v\in \Z^n\setminus\{0\}$.  In particular, 
\begin{align*}
\epsilon \leq \nn{e_n uak} \leq \sqrt{n} \nn{e_n ua} = \sqrt{n}\nn{e_n a} =\sqrt{n}a_n\nn{e_n} = \sqrt{n}a_n.
\end{align*}
But since $a\in A_{2/\sqrt{3}}$, we can also say
\begin{align*}
\epsilon/\sqrt{n} \leq a_n \leq (2/\sqrt{3})a_{n-1}\leq (2/\sqrt{3})^2a_{n-2}\leq\cdots\leq(2/\sqrt{3})^{n-1}a_{1}
\end{align*}
which means that $a_i \geq C\epsilon$ for all $1\leq i \leq n$, where $ C = (\sqrt{3}/2)^{n-1}/\sqrt{n}$.  Moreover, since $\det a = a_1 a_2 \cdots a_n = 1$, we have that
\begin{align*}
a_i = \frac{1}{a_1\cdots a_{i-1}a_{i+1}\cdots a_n} \leq \frac{1}{C^{n-1}\epsilon^{n-1}}
\end{align*}
for any $1\leq i\leq n$.
This implies that for any $1\leq i,j\leq n$, the ratio $a_i/a_j$ can be bounded by
\begin{align*}
\frac{a_i}{a_j} \leq \frac{1}{C^n\epsilon^n}.
\end{align*}
But notice that for an arbitrary matrix $m\in {\rm  Mat}_{n\times n}(\R)$, 
\begin{align*}
|(ama^{-1})_{ij}| = \frac{a_i}{a_j} |m_{ij}| \leq C^{-n}\epsilon^{-n} |m_{ij}|.
\end{align*}
Thus for $a\in A_{2/\sqrt{3}}$, under the max norm on matrices, we have
\begin{align*}
\nn{ama^{-1}}\leq C^{-n} \epsilon^{-n}\nn{m}.
\end{align*}
Furthermore, since $u\in U_{1/2}$, the magnitudes of all entries of $u$ are bounded by $1$.  It is therefore relatively straightforward to see (via matrix multiplication) that $|(umu^{-1})_{ij}|\leq n^2 \max_{i,j} |m_{ij}|$, hence $\nn{umu^{-1}}\leq n^2 \nn{m}$, and the same follows for $k\in K$.
Thus for arbitrary $m\in  {\rm  Mat}_{n\times n}(\R^n)$,
\begin{align*}
\nn{gmg^{-1}} &= \nn{uakmk^{-1}a^{-1}u^{-1}}\\
&\ll \nn{akmk^{-1}a^{-1}}\\
&\ll \epsilon^{-n}\nn{kmk^{-1}}\\
&\ll \epsilon^{-n}\nn{m}
\end{align*}
where all of the above constants depend solely on $n$.  This implies that
\begin{align*}
\nn{{\rm  Ad}(g)}\ll \epsilon^{-n}
\end{align*}
as claimed.
\end{proof}

We may now prove our priginal lemma for the radius of injection.

\begin{proof}[Proof of Lemma \ref{lem:radiusofinjection}]
Let $x\in L_\epsilon$.  By Section \ref{sect:SiegelSet}, we can write $x=\Gamma g$, for some $g\in\Sigma_{1/2,2/\sqrt{3}}$.  Suppose $g_1, g_2\in B^G_r(e)$ and $\pi_x(g_1)=\pi_x(g_2)$, i.e. $\Gamma g g_1 = \Gamma g g_2$.  
Then there exists $\gamma\in\Gamma$ such that $gg_1 = \gamma g g_2$, i.e. $g_1 = g^{-1}\gamma g g_2$.  From this and left-invariance of  the metric, we have that
\begin{align*}
d_G(e,g^{-1}\gamma g)
&\leq d_G(e, g_1) + d_G(g_1,g^{-1}\gamma g)\\
&\leq r + d_G(g\gamma^{-1} g^{-1}g_1,e)\\
&\leq r+ d_G(g\gamma^{-1} g^{-1}g_1,g_2) + d_G(g_2,e)\\
&\leq r+d_G(g_1,g^{-1}\gamma gg_2) + r\\
&= 2r.
\end{align*}

But recall that around every point in $G$ there is a neighborhood on which the metric $d_G$ and the metric derived from any matrix norm are Lipschitz equivalent.  Hence, around the identity, for $r$ less than some fixed value depending only on $n$, we have
\begin{align*}
\nn{e-g^{-1}\gamma g}\ll d_G(e,g^{-1}\gamma g) \ll r
\end{align*}
where $\nn{\cdot}$ is the max norm.  Finally, by Lemma \ref{lem:Adg},
\begin{align*}
\nn{e-\gamma} &=\nn{gg^{-1}(e-\gamma)gg^{-1}}\ll \epsilon^{-n}\nn{g^{-1}(e-\gamma)g}=\epsilon^{-n}\nn{e-g^{-1}\gamma g} \ll r/\epsilon^n.
\end{align*}
Thus for a correctly chosen constant $c_2$, $r=c_2\epsilon^n$ implies that
\begin{align*}
\nn{e-\gamma}<1.
\end{align*}
But since $\gamma\in\Gamma = \slnz$ has integer entries, this can only happen if $\gamma = e$, which implies $g_1=g_2$, so $\pi_x$ is injective on $B^G_r(e)$.
\end{proof}

\section{Properties of the Function $G_d$}\label{app:Gd}

Recall that we defined the generalized Pillai's function $G_d:\N\to\N$ by
\[
G_d (K) := \#\{ {\bf  k} \in \tilde B_K \hspace{2pt}|\hspace{2pt}  k_1 \cdots  k_d \equiv 0 \mod K \}.
\]
We want to prove the following properties of this function.

\begin{customlem}{\ref{lem:Gd}}
 For any integers $K,d\geq 1$, the following hold:
\begin{enumerate}[(i)]
\item (Iterated sum formula) 
\[
G_d(K) = \sum^K_{k_{d-1} = 1}\cdots \sum^K_{k_1 = 1} \gcd(K, k_1\cdots k_{d-1}).
\]
\item (Recursive formula) Let ${\rm  Id}^d(K) = K^d$.  Then
\[
G_{d+1} = {\rm  Id}^d\ast(\phi\cdot G_d).
\]
\item $G_d$ is multiplicative.
\item (Behavior at primes) Let $p$ be a prime.  Then 
\[
G_d(p) = p^d-(p-1)^d.
\]
\item (Dirichlet series bound) For real $x>e$ and $s<d$, 
\[
\sum_{K\leq x} \frac{G_d(K)}{K^s} \ll_{s,d} x^{d-s}(\log x)^{d-1}.
\]
\end{enumerate}
\end{customlem}

To do this, let us first recall a few basic facts from number theory.  For any function $f: \N \to \R$, we have
\begin{align}
\sum^K_{i=1} f(\gcd(K,i)) = \sum_{j|K} f(j)\phi(K/j)\label{eq:Cformula}
\end{align}
where $\phi$ is Euler's totient function, i.e. $\phi(n)$ is the number of positive integers less than $n$ that are relatively prime with $n$. This formula dates back to the work of Ces\`{a}ro and is sometimes refered to as Ces\`{a}ro's formula  (cf. \cite{Cesaro} or \cite{DicksonHist}).

Recall also the definition of Dirichlet convolution: If $f$ and $g$ are functions on the natural numbers, then their convolution is defined by
\begin{align*}
(f\ast g)(K) := \sum_{j|K} f(j)g(K/j).
\end{align*}
So, for example, (\ref{eq:Cformula}) says that $\sum^K_{i=1} f(\gcd(K,i)) = (f\ast \phi)(K)$.  Recall that the convolution of two multiplicative arithmetic functions is again multiplicative. We now have everything we need to complete the proof.

\begin{proof}
%
(\ref{Gdprops1})
To determine an expression for 
\[
G_d (K) := \#\{ {\bf  k} \in \tilde B_K |  k_1 \cdots  k_d \equiv 0 \mod K \}
\]
notice that to specify a point ${\bf k}\in \tilde B_K$ such that $k_1\cdots k_d \equiv 0 \mod K$, we can choose $k_1$ through $k_{d-1}$ independently to be any integers between $1$ and $K$, but then the remaining coordinate $k_d$ must be a multiple of $K/\gcd(K, k_1\cdots k_{d-1})$, that is, the last coordinate must contain all primes in $K$ not contained in any of the previous coordinates.  Since there are $\gcd(K, k_1\cdots k_{d-1})$ multiples of $K/\gcd(K, k_1\cdots k_{d-1})$ less than or equal to $K$, the total number of points counted in this way is given by
\begin{align*}
G_d(K) = \sum^K_{k_{d-1} = 1}\cdots \sum^K_{k_1 = 1} \gcd(K, k_1\cdots k_{d-1}).
\end{align*}

(\ref{Gdprops2})
We will proceed by induction on $d$.  For the base case, we have $G_2(K) = {\rm  Id}\ast \phi = {\rm  Id}\ast (\phi \cdot 1) = {\rm  Id}\ast (\phi \cdot G_1)$, which is a well-known formula for Pillai's arithmetical function, as mentioned above.  Then suppose $G_d = {\rm  Id}^{d-1}\ast (\phi \cdot G_{d-1})$ for $d\geq2$ and consider $G_{d+1}$.

Notice that for any integers $k$,  $n$, and $m$, we can write $\gcd(k, nm) = \gcd(k, n\gcd(k,m))$, that is, we can throw out all the primes in $m$ that are not in $k$.  Furthermore, since $\gcd(k, m)|k$, we can write
\begin{align*}
\gcd(k, n\gcd(k,m)) = \gcd(k, m)\gcd(k/\gcd(k,m), n).
\end{align*}
Hence, $G_{d+1}$ may be written
\begin{align*}
G_{d+1}(K) &= \sum^K_{k_d = 1}\cdots \sum^K_{k_1 = 1} \gcd(K, k_1\cdots k_{d})\\
&= \sum^K_{k_d = 1}\cdots \sum^K_{k_1 = 1} \gcd(K, k_d)\gcd(K/\gcd(K, k_d), k_1\cdots k_{d-1})\\
&= \sum^K_{k_d = 1} \gcd(K, k_d) \left(\sum^K_{k_{d-1} = 1}\cdots \sum^K_{k_1 = 1} \gcd(K/\gcd(K, k_d), k_1\cdots k_{d-1})\right)
\end{align*}
But now notice that the function $\gcd(K/\gcd(K, k_d), k_1\cdots k_{d-1})$ is periodic with period $K/\gcd(K,k_d)$ in each coordinate $k_i$ for $i = 1, \cdots, d-1$.  Thus
\begin{align*}
\sum^K_{k_i = 1} \gcd(K/\gcd(K, k_d), k_1\cdots k_{d-1}) = \gcd(K,k_d)\sum^{K/\gcd(K,k_d)}_{k_i = 1} \gcd(K/\gcd(K, k_d), k_1\cdots k_{d-1})
\end{align*}
for $i = 1, \cdots, d-1$.  Therefore,
\begin{align*}
G_{d+1}(K) &= \sum^K_{k_d=1} \gcd(K, k_d)^d \left(\sum^{K/\gcd(K, k_d)}_{k_{d-1} = 1}\cdots \sum^{K/\gcd(K, k_d)}_{k_1 = 1} \gcd(K/\gcd(K, k_d), k_1\cdots k_{d-1})\right)\\
&= \sum^K_{k_d=1} \gcd(K, k_d)^d G_d(K/\gcd(K, k_d)).
\end{align*}
But by Ces\`{a}ro's formula, this is simply
\begin{align*}
G_{d+1}(K) = \sum_{j|K}  j^d \phi(K/j) G_d(K/j) 
\end{align*}
Finally, we can express this in terms of Dirichlet convolution as
\begin{align*}
G_{d+1}(K) = ( {\rm  Id}^d \ast (\phi \cdot G_d))(K)
\end{align*}
which completes our proof by induction.

(\ref{Gdprops3})
The multiplicativity of $G_d$ for $d\geq 1$ follows immediately from the recursive formula along with the facts that ${\rm  Id}^d$, $\phi$, and $G_1 = 1$ are all multiplicative,
and products and convolutions of multiplicative functions are multiplicative.

(\ref{Gdprops4})
We will again proceed by induction on $d$.  Notice that $G_1(p) = 1 = p-(p-1)$ for all $p$.  Now suppose that for some $d\geq1$, $G_d(p) = p^d - (p-1)^d$ for all primes $p$. 
By the recursive formula proved above, we may write
\begin{align*}
G_{d+1}(p) &= \sum_{j|p}  j^d \phi(p/j) G_d(p/j)\\
&= 1^d\phi(p)G_d(p) +p^d\phi(1) G_d(1)
\end{align*}
since for prime $p$  the sum is only over $j=1,p$.  Then by the induction hypothesis and the facts that $G_d(1) = \phi(1) =1$ and $\phi(p) = p-1$ for any prime $p$, we have that
\begin{align*}
G_{d+1}(p) &= (p-1)(p^d - (p-1)^d) + p^d\\
&= p^{d+1} - (p-1)^{d+1}
\end{align*}
which completes the proof.

(\ref{Gdprops5})
Once again, we proceed by induction on $d$.  Observe that for $d=1$, we have
\begin{align*}
\sum_{K\leq x} \frac{G_1(K)}{K^s} &= \sum_{K\leq x} \frac{1}{K^s}.
\end{align*}
When $0\leq s <1$, we have that $1/K^s$ is decreasing, and $\sum_{K\leq x} 1/K^s \leq 1+\int_1^x 1/t^s dt$.  On the other hand, when $s<0$, we have that $1/K^s$ is increasing, and $\sum_{K\leq x} 1/K^s \leq \int_1^{x+1} 1/t^s dt$.  In either case, we have
\[
\sum_{K\leq x} \frac{G_1(K)}{K^s} \ll_s x^{1-s}
\]
which is the desired bound for $d=1$.

Now suppose that for $d\geq 1$ we have $\sum_{K\leq x} G_d(K)/K^s \ll_{s,d} x^{d-s}(\log x)^{d-1}$ for all $x>e$ and $s<d$.  By the complete multiplicativity of ${\rm  Id}^{-s}$ and the recursive formula for $G_d$, we can write
\[
G_{d+1}(K)/K^s = ({\rm  Id}^{d-s}\ast ({\rm  Id}^{-s}\cdot \phi \cdot G_d))(K).
\]
Also note that a Dirichlet product $(f\ast g)(K) = \sum_{j|K} f(j)g(K/j)$ can be seen as a sum over pairs of positive integers $(n,m)$ whose product is $K$, i.e. 
\[
(f\ast g)(K) = \sum_{\substack{n,m\\nm=K}}f(n)g(m).
\]
Hence, the sum
\[
\sum_{K\leq x} (f\ast g)(K) = \sum_{\substack{n,m\\nm\leq x}} f(n)g(m) = \sum_{n\leq x} f(n)\sum_{m\leq x/n} g(m)
\]
is a sum over pairs of integers whose product is no greater than $x$.  Also notice that $\phi(n)< n$ for any positive integer $n$.  Thus for any $s<d+1$ and $x>e$, we may write
\begin{align*}
\sum_{K\leq x} \frac{G_{d+1}(K)}{K^s} &= \sum_{K\leq x}({\rm  Id}^{d-s}\ast ({\rm  Id}^{-s}\cdot \phi \cdot G_d))(K)\\
&= \sum_{n\leq x}\frac{1}{n^{s-d}} \sum_{m\leq x/n}\frac{\phi(m)G_d(m)}{m^s}\\
&< \sum_{n\leq x}\frac{1}{n^{s-d}} \sum_{m\leq x/n}\frac{G_d(m)}{m^{s-1}}.
\end{align*}
Then since $s<d+1$, we have $s-1<d$.  Also, notice that for $n<x/e$, we have $x/n>e$, so the induction hypothesis applies to sums over $m\leq x/n$ for $n$ in this region.  On the other hand, for $n\geq x/e$, we have $x/n \leq e$, so a sum over $m\leq x/n$ is only a sum over the first two terms, $m=1$ and $m=2$, and can thus be bounded by a constant (depending on $s$ and $d$).  Hence, we may write
\begin{align*}
\sum_{K\leq x} \frac{G_{d+1}(K)}{K^s} &< \sum_{n\leq x}\frac{1}{n^{s-d}} \sum_{m\leq x/n}\frac{G_d(m)}{m^{s-1}}\\
&= \sum_{n< x/e}\frac{1}{n^{s-d}} \sum_{m\leq x/n}\frac{G_d(m)}{m^{s-1}} +\sum_{x/e\leq n\leq x}\frac{1}{n^{s-d}} \sum_{m\leq x/n}\frac{G_d(m)}{m^{s-1}}\\
&\ll_{s,d}  \sum_{n< x/e}\frac{1}{n^{s-d}} (x/n)^{d+1-s}\log(x/n)^{d-1} +\sum_{x/e\leq n\leq x}\frac{1}{n^{s-d}}\\ 
&\ll_{s,d} x^{d+1-s}\sum_{n< x/e}\frac{\log(x/n)^{d-1}}{n} + \sum_{x/e\leq n\leq x}\frac{1}{n^{s-d}}.
\end{align*}
Observe that $\sum_{x/e\leq n\leq x}\frac{1}{n^{s-d}} \ll_{s,d} x^{d+1-s}$ (this can be seen with a calculation similar to that of the base case).  On the other hand, the function $\log(x/t)^{d-1}/t$ is positive and decresing in the region $(1, x/e)$, so we may bound the sum by the first term plus the corresponding integral:
\begin{align*}
\sum_{n< x/e}\frac{\log(x/n)^{d-1}}{n} \leq (\log x)^{d-1} + \int_1^{x/e} \frac{\log(x/t)^{d-1}}{t} dt.
\end{align*}
With the substitution $u=\log(x/t)$, we find that
\[
\int_1^{x/e} \frac{\log(x/t)^{d-1}}{t} dt = \int_1^{\log x}u^{d-1} du = \frac{(\log x)^d -1}{d}.
\]
In total, we have that
\begin{align*}
\sum_{K\leq x} \frac{G_{d+1}(K)}{K^s} &\ll_{s,d} x^{d+1-s} \left(1+(\log x)^{d-1} + (\log x)^d\right)\\
&\ll  x^{d+1-s}(\log x)^d
\end{align*}
since $x>e$, and this completes the proof.
\end{proof}

\end{appendix}


\section*{Acknowledgements}

I would like to thank my advisor, Amir Mohammadi, for giving me the idea to work on this problem and for providing important guidance throughout the process.  I also want to thank Manfred Einsiedler and Hee Oh for helpful discussions about this work.

\bibliographystyle{amsplain}
\bibliography{ResearchDocBib}

\providecommand{\bysame}{\leavevmode\hbox to3em{\hrulefill}\thinspace}
\providecommand{\MR}{\relax\ifhmode\unskip\space\fi MR }
\providecommand{\MRhref}[2]{%
  \href{http://www.ams.org/mathscinet-getitem?mr=#1}{#2}
}
\providecommand{\href}[2]{#2}
\begin{thebibliography}{10}

\bibitem{OEIS}
\emph{A018804}, On-line Encyclopedia of Integer Sequences,
  \url{https://oeis.org/A018804} ({A}ccessed: Feb. 2, 2018).

\bibitem{BekkaMayer}
M.~B. Bekka and M.~Mayer, \emph{Ergodic theory and topological dynamics of
  group actions on homogeneous spaces}, London Mathematical Society Lecture
  Note Series, vol. 269, Cambridge University Press, Cambridge, 2000.
  \MR{1781937}

\bibitem{HeeBenoist}
Y.~Benoist and H.~Oh, \emph{Effective equidistribution of {$S$}-integral points
  on symmetric varieties}, Ann. Inst. Fourier (Grenoble) \textbf{62} (2012),
  no.~5, 1889--1942. \MR{3025156}

\bibitem{BKMNondiv}
V.~Bernik, D.~Y. Kleinbock, and G.~A. Margulis, \emph{Khintchine-type theorems
  on manifolds: the convergence case for standard and multiplicative versions},
  Internat. Math. Res. Notices (2001), no.~9, 453--486. \MR{1829381}

\bibitem{BordGcdGeneral}
O.~Bordell\`es, \emph{Mean values of generalized gcd-sum and lcm-sum
  functions}, J. Integer Seq. \textbf{10} (2007), no.~9, Article 07.9.2, 13.
  \MR{2346091}

\bibitem{BordPillai}
\bysame, \emph{A note on the average order of the gcd-sum function}, J. Integer
  Seq. \textbf{10} (2007), no.~3, Article 07.3.3, 4. \MR{2291947}

\bibitem{Broug1}
K.~A. Broughan, \emph{The gcd-sum function}, J. Integer Seq. \textbf{4} (2001),
  no.~2, Article 01.2.2, 19. \MR{1873400}

\bibitem{Broug2}
\bysame, \emph{The average order of the {D}irichlet series of the gcd-sum
  function}, J. Integer Seq. \textbf{10} (2007), no.~4, Article 07.4.2, 6.
  \MR{2304360}

\bibitem{Burger}
M.~Burger, \emph{Horocycle flow on geometrically finite surfaces}, Duke Math.
  J. \textbf{61} (1990), no.~3, 779--803. \MR{1084459}

\bibitem{Cesaro}
E.~Cesaro, \emph{{\'E}tude moyenne di plus grand commun diviseur de deux
  nombres}, Annali di Matematica Pura ed Applicata (1867-1897) \textbf{13}
  (1885), no.~1, 235--250.

\bibitem{LongGcdGeneral}
J.~Chidambaraswamy and R.~Sitaramachandra~Rao, \emph{Asymptotic results for a
  class of arithmetical functions}, Monatsh. Math. \textbf{99} (1985), no.~1,
  19--27. \MR{778167}

\bibitem{KDabb}
K.~Dabbs, M.~Kelly, and H.~Li, \emph{Effective equidistribution of translates
  of maximal horospherical measures in the space of lattices}, J. Mod. Dyn.
  \textbf{10} (2016), 229--254. \MR{3538863}

\bibitem{Dani78a}
S.~G. Dani, \emph{Invariant measures of horospherical flows on noncompact
  homogeneous spaces}, Invent. Math. \textbf{47} (1978), no.~2, 101--138.
  \MR{0578655}

\bibitem{Dani81}
\bysame, \emph{Invariant measures and minimal sets of horospherical flows},
  Invent. Math. \textbf{64} (1981), no.~2, 357--385. \MR{629475}

\bibitem{Dani84a}
\bysame, \emph{On orbits of unipotent flows on homogeneous spaces}, Ergodic
  Theory Dynam. Systems \textbf{4} (1984), no.~1, 25--34. \MR{758891}

\bibitem{DaniNondiv}
\bysame, \emph{On orbits of unipotent flows on homogeneous spaces. {II}},
  Ergodic Theory Dynam. Systems \textbf{6} (1986), no.~2, 167--182. \MR{857195}

\bibitem{Dani86}
\bysame, \emph{Orbits of horospherical flows}, Duke Math. J. \textbf{53}
  (1986), no.~1, 177--188. \MR{835804}

\bibitem{DaniMarg}
S.~G. Dani and G.~A. Margulis, \emph{Asymptotic behaviour of trajectories of
  unipotent flows on homogeneous spaces}, Proc. Indian Acad. Sci. Math. Sci.
  \textbf{101} (1991), no.~1, 1--17. \MR{1101994}

\bibitem{DicksonHist}
L.~E. Dickson, \emph{History of the theory of numbers. {V}ol. {I}:
  {D}ivisibility and primality}, Chelsea Publishing Co., New York, 1966.
  \MR{0245499}

\bibitem{EMMV}
M.~Einsiedler, G.~A. Margulis, A.~Mohammadi, and A.~Venkatesh, \emph{Effective
  equidistribution and property tau}, (preprint) arXiv:1503.05884 (2015).

\bibitem{EMV}
M.~Einsiedler, G.~A. Margulis, and A.~Venkatesh, \emph{Effective
  equidistribution for closed orbits of semisimple groups on homogeneous
  spaces}, Invent. Math. \textbf{177} (2009), no.~1, 137--212. \MR{2507639}

\bibitem{EW}
M.~Einsiedler and T.~Ward, \emph{Ergodic theory with a view towards number
  theory}, Graduate Texts in Mathematics, vol. 259, Springer-Verlag London,
  Ltd., London, 2011. \MR{2723325}

\bibitem{REWP}
R.~Ellis and W.~Perrizo, \emph{Unique ergodicity of flows on homogeneous
  spaces}, Israel J. Math. \textbf{29} (1978), no.~2-3, 276--284. \MR{0473095}

\bibitem{EMMQuantOpp}
A.~Eskin, G.~A. Margulis, and S.~Mozes, \emph{On a quantitative version of the
  {O}ppenheim conjecture}, Electron. Res. Announc. Amer. Math. Soc. \textbf{1}
  (1995), no.~3, 124--130. \MR{1369644}

\bibitem{QuantOpp}
\bysame, \emph{Upper bounds and asymptotics in a quantitative version of the
  {O}ppenheim conjecture}, Ann. of Math. (2) \textbf{147} (1998), no.~1,
  93--141. \MR{1609447}

\bibitem{FF}
L.~Flaminio and G.~Forni, \emph{Invariant distributions and time averages for
  horocycle flows}, Duke Math. J. \textbf{119} (2003), no.~3, 465--526.
  \MR{2003124}

\bibitem{Furst}
H.~Furstenberg, \emph{The unique ergodicity of the horocycle flow}, Recent
  advances in topological dynamics ({P}roc. {C}onf., {Y}ale {U}niv., {N}ew
  {H}aven, {C}onn., 1972; in honor of {G}ustav {A}rnold {H}edlund), Springer,
  Berlin, 1973, pp.~95--115. \MR{0393339}

\bibitem{Conjectures}
A.~Gorodnik, \emph{Open problems in dynamics and related fields}, J. Mod. Dyn.
  \textbf{1} (2007), no.~1, 1--35. \MR{2261070}

\bibitem{GMO}
A.~Gorodnik, F.~Maucourant, and H.~Oh, \emph{Manin's and {P}eyre's conjectures
  on rational points and adelic mixing}, Ann. Sci. \'Ec. Norm. Sup\'er. (4)
  \textbf{41} (2008), no.~3, 383--435. \MR{2482443}

\bibitem{GreenTao1}
B.~Green and T.~Tao, \emph{The quantitative behaviour of polynomial orbits on
  nilmanifolds}, Ann. of Math. (2) \textbf{175} (2012), no.~2, 465--540.
  \MR{2877065}

\bibitem{CombSieve}
H.~Halberstam and H.-E. Richert, \emph{Sieve methods}, Academic Press [A
  subsidiary of Harcourt Brace Jovanovich, Publishers], London-New York, 1974,
  London Mathematical Society Monographs, No. 4. \MR{0424730}

\bibitem{Heegner2}
G.~Harcos and P.~Michel, \emph{The subconvexity problem for {R}ankin-{S}elberg
  {$L$}-functions and equidistribution of {H}eegner points. {II}}, Invent.
  Math. \textbf{163} (2006), no.~3, 581--655. \MR{2207235}

\bibitem{HaukGcdGeneral}
P.~Haukkanen, \emph{On a gcd-sum function}, Aequationes Math. \textbf{76}
  (2008), no.~1-2, 168--178. \MR{2443468}

\bibitem{Hedl}
G.~A. Hedlund, \emph{Fuchsian groups and transitive horocycles}, Duke Math. J.
  \textbf{2} (1936), no.~3, 530--542. \MR{1545946}

\bibitem{KMBddOrbits}
D.~Y. Kleinbock and G.~A. Margulis, \emph{Bounded orbits of nonquasiunipotent
  flows on homogeneous spaces}, Sina\u\i 's {M}oscow {S}eminar on {D}ynamical
  {S}ystems, Amer. Math. Soc. Transl. Ser. 2, vol. 171, Amer. Math. Soc.,
  Providence, RI, 1996, pp.~141--172. \MR{1359098}

\bibitem{KMNondiv}
\bysame, \emph{Flows on homogeneous spaces and {D}iophantine approximation on
  manifolds}, Ann. of Math. (2) \textbf{148} (1998), no.~1, 339--360.
  \MR{1652916}

\bibitem{KMLogLaws}
\bysame, \emph{Logarithm laws for flows on homogeneous spaces}, Invent. Math.
  \textbf{138} (1999), no.~3, 451--494. \MR{1719827}

\bibitem{KMExpanding}
\bysame, \emph{On effective equidistribution of expanding translates of certain
  orbits in the space of lattices}, Number theory, analysis and geometry,
  Springer, New York, 2012, pp.~385--396. \MR{2867926}

\bibitem{Knapp}
A.~W. Knapp, \emph{Lie groups beyond an introduction}, Progress in Mathematics,
  vol. 140, Birkh\"auser Boston, Inc., Boston, MA, 1996. \MR{1399083}

\bibitem{LeeOh}
M.~Lee and H.~Oh, \emph{Effective equidistribution of closed horocycles for
  geometrically finite surfaces}, (preprint) arXiv:1202.0848 (2012).

\bibitem{Marg71}
G.~A. Margulis, \emph{On the action of unipotent groups in a lattice space},
  Mat. Sb. (N.S.) \textbf{86(128)} (1971), no.~4(12), 552--556. \MR{0291352}

\bibitem{MargulisNondiv}
\bysame, \emph{On the action of unipotent groups in the space of lattices}, Lie
  groups and their representations ({P}roc. {S}ummer {S}chool, {B}olyai,
  {J}\'anos {M}ath. {S}oc., {B}udapest, 1971), Halsted, New York, 1975,
  pp.~365--370. \MR{0470140}

\bibitem{Opp}
\bysame, \emph{Formes quadratriques ind\'efinies et flots unipotents sur les
  espaces homog\`enes}, C. R. Acad. Sci. Paris S\'er. I Math. \textbf{304}
  (1987), no.~10, 249--253. \MR{882782}

\bibitem{MOpp}
\bysame, \emph{Discrete subgroups and ergodic theory}, Number theory, trace
  formulas and discrete groups ({S}ymposium in honor of {A}tle {S}elberg,
  {O}slo, {N}orway, {J}uly 14-21, 1987), Academic Press, Boston, MA, 1989,
  pp.~377--398. \MR{993328}

\bibitem{Mthesis}
\bysame, \emph{On some aspects of the theory of {A}nosov systems}, Springer
  Monographs in Mathematics, Springer-Verlag, Berlin, 2004, Translated from
  Russian by V. V. Szulikowska. \MR{2035655}

\bibitem{MT}
G.~A. Margulis and G.~M. Tomanov, \emph{Invariant measures for actions of
  unipotent groups over local fields on homogeneous spaces}, Invent. Math.
  \textbf{116} (1994), no.~1-3, 347--392. \MR{1253197}

\bibitem{Heegner1}
P.~Michel, \emph{The subconvexity problem for {R}ankin-{S}elberg
  {$L$}-functions and equidistribution of {H}eegner points}, Ann. of Math. (2)
  \textbf{160} (2004), no.~1, 185--236. \MR{2119720}

\bibitem{SarnackNevo}
A.~Nevo and P.~Sarnak, \emph{Prime and almost prime integral points on
  principal homogeneous spaces}, Acta Math. \textbf{205} (2010), no.~2,
  361--402. \MR{2746350}

\bibitem{Pillai}
S.~S. Pillai, \emph{On an arithmetic function}, J. of the Annamalai Univ.
  \textbf{2}.

\bibitem{Rag}
M.~S. Raghunathan, \emph{Discrete subgroups of {L}ie groups}, Springer-Verlag,
  New York-Heidelberg, 1972, Ergebnisse der Mathematik und ihrer Grenzgebiete,
  Band 68. \MR{0507234}

\bibitem{Rat3}
M.~Ratner, \emph{On {R}aghunathan's measure conjecture}, Ann. of Math. (2)
  \textbf{134} (1991), no.~3, 545--607. \MR{1135878}

\bibitem{RatSummary}
\bysame, \emph{Invariant measures and orbit closures for unipotent actions on
  homogeneous spaces}, Geom. Funct. Anal. \textbf{4} (1994), no.~2, 236--257.
  \MR{1262705}

\bibitem{SarnackAsymp}
P.~Sarnak, \emph{Asymptotic behavior of periodic orbits of the horocycle flow
  and {E}isenstein series}, Comm. Pure Appl. Math. \textbf{34} (1981), no.~6,
  719--739. \MR{634284}

\bibitem{SarnackUbis}
P.~Sarnak and A.~Ubis, \emph{The horocycle flow at prime times}, J. Math. Pures
  Appl. (9) \textbf{103} (2015), no.~2, 575--618. \MR{3298371}

\bibitem{ShahConj}
N.~A. Shah, \emph{Limit distributions of polynomial trajectories on homogeneous
  spaces}, Duke Math. J. \textbf{75} (1994), no.~3, 711--732. \MR{1291701}

\bibitem{StromClosed}
A.~Str\"ombergsson, \emph{On the uniform equidistribution of long closed
  horocycles}, Duke Math. J. \textbf{123} (2004), no.~3, 507--547. \MR{2068968}

\bibitem{StromDev}
\bysame, \emph{On the deviation of ergodic averages for horocycle flows}, J.
  Mod. Dyn. \textbf{7} (2013), no.~2, 291--328. \MR{3106715}

\bibitem{StromRat}
\bysame, \emph{An effective {R}atner equidistribution result for {${\rm
  SL}(2,\mathbb{R})\ltimes\mathbb{R}^2$}}, Duke Math. J. \textbf{164} (2015),
  no.~5, 843--902. \MR{3332893}

\bibitem{Tanigawa}
Y.~Tanigawa and W.~Zhai, \emph{On the gcd-sum function}, J. Integer Seq.
  \textbf{11} (2008), no.~2, Article 08.2.3, 11. \MR{2413094}

\bibitem{TothGcdGeneral}
L.~T\'oth, \emph{A generalization of {P}illai's arithmetical function involving
  regular convolutions}, Proceedings of the 13th {C}zech and {S}lovak
  {I}nternational {C}onference on {N}umber {T}heory ({O}stravice, 1997),
  vol.~6, University of Ostrava, 1998, pp.~203--217. \MR{1828135}

\bibitem{Toth}
\bysame, \emph{A survey of gcd-sum functions}, J. Integer Seq. \textbf{13}
  (2010), no.~8, Article 10.8.1, 23. \MR{2718232}

\bibitem{TothSimilar}
\bysame, \emph{Another generalization of the gcd-sum function}, Arab. J. Math.
  (Springer) \textbf{2} (2013), no.~3, 313--320. \MR{3079887}

\bibitem{Veech}
W.~A. Veech, \emph{Unique ergodicity of horospherical flows}, Amer. J. Math.
  \textbf{99} (1977), no.~4, 827--859. \MR{0447476}

\bibitem{VenkSparse}
A.~Venkatesh, \emph{Sparse equidistribution problems, period bounds and
  subconvexity}, Ann. of Math. (2) \textbf{172} (2010), no.~2, 989--1094.
  \MR{2680486}

\end{thebibliography}

\end{document}